\documentclass[11pt]{article}
\usepackage{amsthm}
\usepackage{pxfonts}
\usepackage{yfonts}
\usepackage{dsfont}
\usepackage{graphicx}
\usepackage{relsize}
\usepackage{color}

\def\bel{\begin{equation}\label}
\def\eeq{\end{equation}}
\def\ds{\displaystyle}
\def\endproof{\hphantom{MM}
\hfill\llap{$\square$}\goodbreak}

\def\mt{\longrightarrow}
\def\v{\vskip 1em}

\def\ve{\varepsilon}
\def\vex{\epsilon}
\def\R{\mathbb R}
\def\Z{\mathbb Z}
\def\C{\mathfrak{C}}
\def\Cx{\mathds C}

\def\Re{{\bf Re}}
\def\Im {{\bf Im}}

\def\Q{{\bf Q}}

\def\A{{\bf A}}

\def\L{{\bf L}}
\def\U{{\bf U}}
\def\V{{\bf V}}
\def\T{{\bf T}}

\def\i{{\bf i}}

\def\Hat{\widehat}

\def\vol{{\bf vol}}
\def\sign{{\bf sign}}

\def\I{{\bf I}}

\def\M{{\bf M}}
\def\G{{\bf G}}

\def\Cup{{\bigcup}}
\def\Cap{{\bigcap}}

\def\alpha{\alphaup}
\def\beta{\betaup}
\def\gamma{\gammaup}
\def\delta{\deltaup}
\def\theta{\thetaup}
\def\xi{{\xiup}}
\def\eta{{\etaup}}
\def\tau{{\tauup}}
\def\rho{{\rhoup}}
\def\phi{{\phiup}}
\def\psi{{\psiup}}
\def\lambda{{\lambdaup}}
\def\omega{\omegaup}
\def\varphi{{\varphiup}}
\def\gamma{{\gammaup}}

\def\j{\jmath}

\newtheorem{lemma}{Lemma}[section]

\newtheorem{remark}{Remark}[section]

\begin{document}
\[\hbox{\LARGE{\bf Fractional integrals associated with Zygmund dilations}}\]

\[\hbox{Zipeng Wang}\]
\begin{abstract}
We study a family of fractional integral operators defined in $\R^3$ whose kernels are distributions associated with Zygmund dilations: $(x_1,x_2,x_3)\mt(\delta_1 x_1, \delta_2 x_2, \delta_1\delta_2 x_3)$ for  $\delta_1,\delta_2>0$ having singularity on every coordinate subspace. As a result, we obtain a Hardy-Littlewood-Sobolev type inequality.
\end{abstract}

\section{Introduction}
\setcounter{equation}{0}
In 1928, Hardy and Littlewood \cite{Hardy-Littlewood} introduced a family of convolution operators
\bel{i_alpha}
\begin{array}{cc}\ds
I_\alpha f(x)~=~\int_{\R} f(y)\left|x-y\right|^{\alpha-1} dy,\qquad 0<\alpha<1.
 \end{array}
\eeq
$\diamond$ {\small Throughout, $\C>0$ is regarded as a generic constant whose value depends on the sub-indices.}

{\bf Theorem A: Hardy and Littlewood, 1928} ~~~{\it Let $I_\alpha$ defined in (\ref{i_alpha}). We have
\bel{ineq 1} \left\| I_\alpha f\right\|_{\L^q(\R)}~\leq~\C_{p~q}~\left\| f\right\|_{\L^p(\R)},\qquad 1<p<q<\infty
\eeq
if and only if 
\bel{formula 1}
\alpha~=~{1\over p}-{1\over q}.
\eeq
}\\
Ten years later, this result has been extended to higher dimensional spaces by Sobolev \cite{Sobolev}. Today, it bears  the name of  Hardy-Littlewood-Sobolev inequality.

The equation in  (\ref{formula 1}) is called the homogeneity condition of (\ref{ineq 1}). As an necessity, (\ref{formula 1}) guarantees that the norm inequality in (\ref{ineq 1}) is invariant by changing dilations: $I_\alpha f(x)\mt I_\alpha f(\delta x)$ and $f(x)\mt f(\delta x)$ for $\delta>0$.

{\bf Theorem A} can be easily extended to the multi-parameter setting. For instance, define
\bel{I 123}
I_{\alpha_1\alpha_2\alpha_3}f(x)~=~\int_{\R^3} f(y)\prod_{i=1}^3 \left|x_i-y_i\right|^{\alpha_i-1}dy,\qquad 0<\alpha_i<1,~~ i=1,2,3.
\eeq
Observe that the kernel of $I_{\alpha_1\alpha_2\alpha_3}$ has singularity at every $x_i=0, i=1,2,3$.

{\bf Theorem B:} ~~~{\it Let $I_{\alpha_1\alpha_2\alpha_3}$ defined in (\ref{I 123}). We have
\bel{ineq 2} \left\| I_{\alpha_1\alpha_2\alpha_3} f\right\|_{\L^q(\R^3)}~\leq~\C_{p~q}~\left\| f\right\|_{\L^p(\R^3)},\qquad 1<p<q<\infty
\eeq
if and only if 
\bel{formula 2}
\alpha_1~=~\alpha_2~=~\alpha_3~=~{1\over p}-{1\over q}.
\eeq
}\\
The required homogeneity condition in (\ref{formula 2}) can be shown by carrying out a $3$-parameter changing dilations: $I_{\alpha_1\alpha_2\alpha_3}f( x)\mt I_{\alpha_1\alpha_2\alpha_3}f(\delta_1 x_1, \delta_2 x_2, \delta_3 x_3)$
and $f(x)\mt f(\delta_1 x_1, \delta_2 x_2,\delta_3 x_3)$ for $\delta_1, \delta_2, \delta_3>0$ inside (\ref{ineq 2}).
Conversely, we prove the norm inequality in (\ref{ineq 2}) by using a familiar iteration argument. \footnote{We apply {\bf Theorem A} in each coordinate subspace together with using Minkowski integral inequality.}

In this paper, we investigate a new type of fractional integral operators whose kernels are distributions associated with Zygmund dilations.  This is
 a group of dilations in $\R^3$ lying between the standard $1$-parameter dilations and the $3$-parameter dilations discussed  above. 
Namely, we assert $(x_1,x_2,x_3)\mt(\delta_1 x_1, \delta_2 x_2, \delta_1\delta_2 x_3)$ for $\delta_1, \delta_2>0$.

\subsection{Maximal functions, Singular integrals and Zygmund dilations }
Let $\hbox{\bf R}=Q_1\times Q_2\times Q_3$ denote a rectangle in $\R^3$ where each $Q_i, i=1,2,3$ is an open interval. The maximal function commute with Zygmund dilations is defined as
\bel{M_z}
\M_\zeta f(x)~=~\sup_{x\in\hbox{\bf R}:~|Q_3|=|Q_1||Q_2|}{1\over |\hbox{\bf R}|}\int_{\hbox{\bf R}} |f(y)|dy.
\eeq
Stein was the first to link the properties of $\M_\zeta$ to the boundary value problem of Poisson integrals in symmetric spaces, such as Siegel's upper half space.   We refer to the survey paper of Fefferman \cite{Fefferman} for more discussions on $\M_\zeta$.

Another concrete example is given by Nagel and Wainger \cite{Nagel-Wainger} considering convolutions with the kernel $\Omega$ that agrees on the function
\bel{Nagel-Wainger kernel}
\Omega(x)~=~\sign(x_1x_2) |x_1|^{-1}|x_2|^{-1}|x_3|^{-1} \Bigg[ {|x_1||x_2|\over |x_3|}+{|x_3|\over |x_1||x_2|}\Bigg]^{-1}
\eeq
away from the subspace $x_1=x_3=0$ or $x_2=x_3=0$. The $\L^2$-boundedness of $f\ast\Omega$ is obtained as a special case among the results in \cite{Nagel-Wainger}. 

Our next example is related  to the usual representation of the Heisenberg group in $\R^3$ with group law
$(x_1,x_2,x_3)\odot (y_1,y_2,y_3)=\Big[x_1+y_1, x_2+y_2, x_3+y_3+\mu(x_1y_2-y_1x_2)\Big]$, $ \mu\neq0$.

The regarding Cauchy-Szego's kernel is given by
\bel{Cauchy-Szego}
\mathfrak{S}(x)~=~\left[{1\over x_3+\i (x_1^2+x_2^2)}\right]^2,\qquad x\neq0.
\eeq
Convolutions with $\mathfrak{S}$ in the Heisenberg group represented in $\R^3$ are defined as
\bel{T}
\T f(x)~=~\int_{\R^3} f(y) \mathfrak{S}(x\odot y^{-1})dy.
\eeq
$\T$ defined in (\ref{T}) is commute with Zygmund dilations whenever $\delta_1=\delta_2$ occurs, i.e: $(x_1,x_2,x_3)\mt (\delta x_1, \delta x_2, \delta^2 x_3)$, $\delta>0$. Furthermore, it is 
bounded on $\L^p(\R^3)$ for $1<p<\infty$. More background can be found in chapter XII of Stein \cite{Stein}. 
See   M\"{u}ller, Ricci and Stein \cite{M.R.S} concerning Fourier multipliers of Marcinkiewicz type in this setting.

In 1992, Ricci and Stein \cite{Ricci-Stein} introduced a general class of singular integral operators  which can be  characterized by their kernels or equivalently by  the  regarding Fourier multipliers. An $\L^p$-theorem  has been established. Later, these singular integral operators 
are refined by Fefferman and Pipher  \cite{Fefferman-Pipher} for the weighted analogous $\L^p$-estimates.

More recently, a larger family of singular integrals associated with Zygmund dilations is invented by Han et al \cite{Han-Li-Lin-Tan}. In compare to the settings of the previous works in \cite{Ricci-Stein} and \cite{Fefferman-Pipher}, the regularity condition assigned on the kernels  has been reduced to involve only a minimal H\"{o}lder-continuity type estimates. The $\L^p$-boundedness of singular integral operators under consideration is also concluded.

The latest update in this direction refers to the beautiful paper by  Hyt\"{o}nen et al \cite{Hytonen-}: A characterization is found between the weighted $\L^p$-norm inequalities for a certain class of singular integrals defined in \cite{Han-Li-Lin-Tan} and the corresponding Muckenhoupt $\A_p$-class satisfying Zygmund dilations. Novel examples are provided to show the optimality of  this special $\A_p$-class $w.r.t$ the regularity and cancellation conditions carried by the kernels.

\subsection{Formulation on the main result}
Besides the substantial development for singular integrals, the area of fractional integrals associated with Zygmund dilations remains largely open.  
Motivated by the explicit examples shown in (\ref{M_z}), (\ref{Nagel-Wainger kernel}) and (\ref{Cauchy-Szego})-(\ref{T}), we consider
\bel{I_alpha123}
\I_{\alpha_1\alpha_2\alpha_3} f(x)~=~\int_{\R^3} f(y) \V^{\alpha_1\alpha_2\alpha_3}(x-y)dy
\eeq
 where
\bel{Kernel}
\V^{\alpha_1\alpha_2\alpha_3}(x)~=~|x_1|^{\alpha_1-1}|x_2|^{\alpha_2-1}|x_3|^{\alpha_3-1} \Bigg[ {|x_1||x_2|\over |x_3|}+{|x_3|\over |x_1||x_2|}\Bigg]^{-1}
\eeq
for $\alpha_i<1$, $i=1,2,3$ whenever $\I_{\alpha_1\alpha_2\alpha_3}$ is well defined. 

Our aim is to find a characterization between the norm inequality
\bel{norm ineq}
\left\| \I_{\alpha_1\alpha_2\alpha_3} f\right\|_{\L^q(\R^3)}~\leq~\C~\left\| f\right\|_{\L^p(\R^3)},\qquad 1<p<q<\infty
\eeq
and the necessary constraints consisting of $\alpha_1, \alpha_2,\alpha_3, p,q$. 
\v

Suppose ${1\over 2}\leq|x_3|<2$ and $|x_1|\leq1, |x_2|\leq1$. We find $\V^{\alpha_1\alpha_2\alpha_3}(x)\approx |x_1|^{\alpha_1}|x_2|^{\alpha_2}$. This implies $\alpha_1>-1, \alpha_2>-1$ due to the essential local integrability of $\V^{\alpha_1\alpha_2\alpha_2}(x)$. On the other hand, suppose ${1\over 2}\leq|x_1|<2$,  ${1\over 2}\leq|x_2|<2$ and $|x_2|\leq1$. We find $\V^{\alpha_1\alpha_2\alpha_3}(x)\approx |x_3|^{\alpha_3}$. Hence that 
 $\alpha_3>-1$ is an necessity.

Next, consider ${1\over 2}\leq|x_1|<2$ and ${1\over 2}|x_2|<|x_3|<2|x_2|$.
We find 
\bel{V 123 est size}
\V^{\alpha_1\alpha_2\alpha_3}(x)~\approx~ |x_2|^{\alpha_2-1}|x_3|^{\alpha_3-1}~\approx~\left[{1\over |x_2|+|x_3|}\right]^{2-\alpha_2-\alpha_3}.
\eeq
This further implies $\alpha_2+\alpha_3>0$ as required for the local integrability of $\V^{\alpha_1\alpha_2\alpha_2}(x)$. Moreover, we also have $\alpha_1+\alpha_3>0$ for symmetry reason.

In summary of the above, we essentially need
\bel{alpha123}
-1<\alpha_i<1,~~ i=1,2,3\qquad\hbox{and}\qquad \alpha_1+\alpha_3>0, \qquad \alpha_2+\alpha_3>0
\eeq
to make $\I_{\alpha_1\alpha_2\alpha_2}$ well defined in (\ref{I_alpha123}).

By changing dilations $\I_{\alpha_1\alpha_2\alpha_2}f(x)\mt \I_{\alpha_1\alpha_2\alpha_2}(\delta_1 x_1, \delta_2 x_2, \delta_1\delta_2 x_3)$ and $f(x)\mt f(\delta_1 x_1, \delta_2 x_2, \delta_1\delta_2 x_3)$
inside (\ref{norm ineq}), the norm inequality  implies 
\bel{12 constraints}
{\alpha_1+\alpha_3\over 2}~=~{1\over p}-{1\over q},\qquad  {\alpha_2+\alpha_3\over 2}~=~{1\over p}-{1\over q}
\eeq
simultaneously. Therefore, we must have $\alpha_1=\alpha_2$. 

Now, we can redefine our fractional integral operator associated with Zygmund dilations.
Let $-1<\alpha,\beta<1$ and $\alpha+\beta>0$. Consider
\bel{I_alpha beta}
\I_{\alpha\beta} f(x)~=~\int_{\R^3} f(y) \V^{\alpha\beta}(x-y)dy
\eeq
where
\bel{V}
\V^{\alpha\beta}(x)~=~|x_1|^{\alpha-1}|x_2|^{\alpha-1}|x_3|^{\beta-1} \Bigg[ {|x_1||x_2|\over |x_3|}+{|x_3|\over |x_1||x_2|}\Bigg]^{-1} 
\eeq
for $x_i\neq0, i=1,2,3$. 
Our main result is stated in below.
\v
{\bf Theorem One}~~{\it Let $\I_{\alpha\beta}$ defined in (\ref{I_alpha beta})-(\ref{V}) for $-1<\alpha,\beta<1$ and $\alpha+\beta>0$. We have
\bel{Result One}
\left\| \I_{\alpha\beta} f\right\|_{\L^q(\R^3)}~\leq~\C_{\alpha~\beta~p~q}~\left\| f\right\|_{\L^p(\R^3)},\qquad 1<p<q<\infty
\eeq
if and only if
\bel{Formula}
{\alpha+\beta\over 2}~=~{1\over p}-{1\over q}.
\eeq}

\subsection{Sketch on the proof of Theorem One}
In the next section, we develop a new framework by asserting $\I_{\alpha\beta}f=\sum_{\ell\in\Z} \Delta_\ell \I_{\alpha\beta}f$ for which
\[ \Delta_\ell \I_{\alpha\beta} f(x)~=~\int_{\Gamma_\ell(x)} f(y) \V^{\alpha\beta}(x-y)dy,\qquad \Cup_{\ell\in\Z} \Gamma_\ell(x)=\R^3.
\]
The projection of  $\Gamma_\ell(x)$ in the $(x_1,x_2)$-subspace is a collection of dyadic rectangles having a same eccentricity depending on $\ell\in\Z$ with side length comparable to the distance from $(x_1,x_2)\in\R^2$. We shall see that every  $\Delta_\ell \I_{\alpha\beta}$ satisfies the desired $\L^p\mt\L^q$-norm inequality. Furthermore, they enjoy a certain almost orthogonality property, stated as {\bf Proposition One}. 

In section 3, we derive a point-wise estimate to dominate each $\Delta_\ell \I_{\alpha\beta}$ by using the strong maximal function, as a multi-parameter analogous of Hedberg \cite{Hedberg}. 
Section 4 is devoted to some implications on {\bf Proposition One}. 
We finish the proof in section 5.
\begin{remark}
Because $\I_{\alpha\beta}$ is positively definite, we will assume $f\ge0$ in the remaining paper.
\end{remark}

\section{Dyadic decomposition in $\R\times\R\times\R$}
\setcounter{equation}{0}
Let $\ell, j,k\in\Z$. Define
\bel{Gamma_ljk}
\Gamma_{\ell j k}(x)~=~\left\{ y\in\R^3\colon \left.\begin{array}{cc}\ds2^j\leq|x_1-y_1|<2^{j+1},~2^{j-\ell}\leq|x_2-y_2|<2^{j+1-\ell},
\\\\ \ds
2^{j+(j-\ell)-k}\leq|x_3-y_3|<2^{j+(j-\ell)+1-k}
\end{array}\right.\right\}
\eeq
and
\bel{Gamma_lj}
\Gamma_\ell(x)~=~\Cup_{j\in\Z}~ \Gamma_{\ell j}(x),\qquad \Gamma_{\ell j}(x)~=~\Cup_{k\in\Z}~ \Gamma_{\ell jk}(x).
\eeq
\begin{figure}[h]
\centering
\includegraphics[scale=0.52]{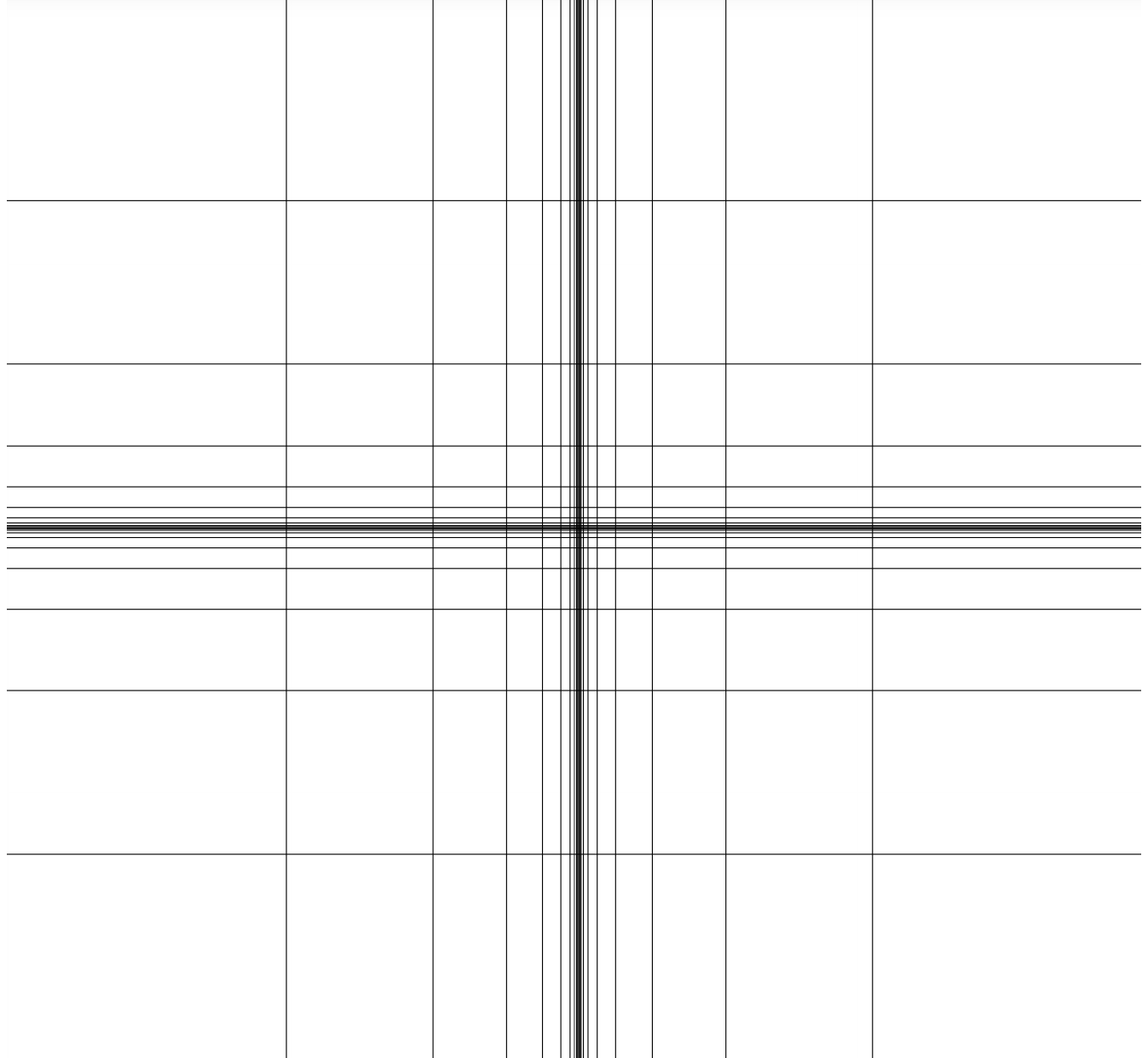}
\caption{Projections of $\Gamma_{\ell j}(x)$, $\ell, j\in\Z$ in $\R^2$}
\end{figure}

$\Gamma_{\ell j}(x)$ is a rectangular cylinder in $\R^3$ which is independent from $x_3$, i.e: $\Gamma_{\ell j}(x)=\Gamma_{\ell j}(x_1,x_2)$. The projection of $\Gamma_{\ell j}(x)$ in the $(x_1,x_2)$-subspace is a dyadic rectangle denoted by
\bel{Gamma project}
\begin{array}{cc}\ds
\Lambda_{\ell j}(x_1,x_2)~=~\Lambda^1_{\ell j}(x_1)\Lambda^2_{\ell j}(x_2),
\\\\ \ds
\Lambda^1_{\ell j}(x_1)~=~\left\{ y_1\in\R\colon 2^{j}\leq|x_1-y_1|<2^{j+1}\right\},\qquad \Lambda^2_{\ell j}(x_2)~=~\left\{ y_2\in\R\colon 2^{j-\ell}\leq|x_2-y_2|<2^{j-\ell+1}\right\}.
\end{array}
\eeq
The projection of $\Gamma_\ell(x)$ in the $(x_1,x_2)$-subspace is a collection of dyadic rectangles having the same eccentricity. Geometrically, it can be interpreted as a discrete version of "dyadic cone" vertex on $(x_1, x_2)$.

We respectively define
\bel{Delta_l I}
\begin{array}{ccc}\ds
\Delta_\ell \I_{\alpha\beta} f(x)~=~\int_{\Gamma_\ell(x)} f(y) \V^{\alpha\beta}(x-y)dy,
\\\\ \ds
\Delta_{\ell j} \I_{\alpha\beta} f(x)~=~\int_{\Gamma_{\ell j}(x)} f(y) \V^{\alpha\beta}(x-y)dy,
\qquad
\Delta_{\ell jk} \I_{\alpha\beta} f(x)~=~\int_{\Gamma_{\ell jk}(x)} f(y) \V^{\alpha\beta}(x-y)dy.
\end{array}
\eeq
\v
{\bf Proposition One}:~~{\it Let ${\alpha+\beta\over 2}={1\over p}-{1\over q},1<p<q<\infty$. Suppose  $q\in\Z$ satisfying $(q-2)\left[{\alpha+\beta\over 2}\right]\ge1$. We have
\bel{Lemma}
\int_{\R^3} \sum_{\ell\in\Z} \Delta_\ell \I_{\alpha\beta} f(x) \Big(\Delta_{\ell-h}\I_{\alpha\beta} f\Big)^{q-1}(x)dx~\leq~\C_{\alpha~\beta~q} ~2^{-\ve |h|}~\left\| f\right\|_{\L^p(\R^3)}^q
\eeq
for some $\ve=\ve(\alpha,\beta,q)>0$.}

Let $h_m\in\Z$ for $m=1,2,\ldots,q-1$. By applying Tonelli's theorem, we write
\bel{I Sum}
\begin{array}{lr}\ds
\int_{\R^3} \Big(\I_{\alpha\beta} f\Big)^q(x)dx
~=~\int_{\R^3} \sum_{h_m\in\Z, ~m=1,2,\ldots,q-1} \sum_{\ell\in\Z} \Delta_\ell \I_{\alpha\beta} f(x) \prod_{m=1}^{q-1}  \Delta_{\ell-h_m} \I_{\alpha\beta} f(x) dx
\\\\ \ds~~~~~~~~~~~~~~~~~~~~~~~~~~~~
~=~\sum_{h_m\in\Z,~ m=1,2,\ldots,q-1}\int_{\R^3}  \sum_{\ell\in\Z} \Delta_\ell \I_{\alpha\beta} f(x) \prod_{m=1}^{q-1}  \Delta_{\ell-h_m} \I_{\alpha\beta} f(x) dx.
\end{array}
\eeq
By using H\"{o}lder inequality twice, we find
\bel{Holder 2 times}
\begin{array}{lr}\ds
\int_{\R^3}  \sum_{\ell\in\Z} \Delta_\ell \I_{\alpha\beta} f(x) \prod_{m=1}^{q-1}  \Delta_{\ell-h_m} \I_{\alpha\beta} f(x) dx~\leq~\int_{\R^3} \prod_{m=1}^{q-1} \left\{\sum_{\ell\in\Z} \Delta_\ell \I_{\alpha\beta} f(x) \Big(\Delta_{\ell-h_m}\I_{\alpha\beta} f\Big)^{q-1}(x)\right\}^{1\over q-1} dx
\\\\ \ds~~~~~~~~~~~~~~~~~~~~~~~~~~~~~~~~~~~~~~~~~~~~~~~~~~~~~~~~~~~~~~
~\leq~\prod_{m=1}^{q-1} \left\{\int_{\R^3} \sum_{\ell\in\Z} \Delta_\ell \I_{\alpha\beta} f(x) \Big(\Delta_{\ell-h_m}\I_{\alpha\beta} f\Big)^{q-1}(x) dx\right\}^{1\over q-1}.
\end{array}
\eeq
From (\ref{I Sum})-(\ref{Holder 2 times}), by applying {\bf Proposition One}, we obtain the norm inequality in (\ref{Result One}) for $q\in\Z$ and $(q-2)\left[{\alpha+\beta\over 2}\right]\ge1$.
Note that $\I_{\alpha\beta}$ defined in (\ref{I_alpha beta}) is  self-adjoint. Therefore, we have $\I_{\alpha\beta}\colon \L^p\left(\R^3\right)\mt\L^q\left(\R^3\right)\Longleftrightarrow
\I_{\alpha\beta}\colon \L^{q\over q-1}\left(\R^3\right)\mt\L^{p\over p-1}\left(\R^3\right)
$. Let
\bel{formula12}
\begin{array}{rl}\ds
{\alpha+\beta\over 2}~=~{1\over p}-{1\over q}~=~{1\over p_1}-{1\over q_1}~=~{1\over p_2}-{1\over q_2},
\qquad
 1~<~p_i~<~q_i~<~\infty,\qquad i~=~1,2.
 \end{array}
\eeq
We choose $q_1,\left({p_2\over p_2-1}\right)\in\Z$  satisfying $(q_1-2)\left[{\alpha+\beta\over 2}\right]\ge1$, $q_1>q$ and $\left({p_2\over p_2-1}-2\right)\left[{\alpha+\beta\over 2}\right]\ge1$, $p_2<p$.
There exists a $0< t< 1$ such that 
\bel{t,pq}
{1\over p}~=~{1-t\over p_1}+{t\over p_2},\qquad{1\over q}~=~{1-t\over q_1}+{t\over q_2}.
\eeq
From the above estimates, we simultaneously have 
\bel{Endpoints Est}
\left\|\I_{\alpha\beta} f\right\|_{\L^{q_1}\left(\R^3\right)}~\leq~\C_{\alpha~\beta~q}~\left\|f\right\|_{\L^{p_1}\left(\R^3\right)},\qquad \left\|\I_{\alpha\beta} f\right\|_{\L^{q_2}\left(\R^3\right)}~\leq~\C_{\alpha~\beta~p}~\left\|f\right\|_{\L^{p_2}\left(\R^3\right)}.
\eeq  
By applying Riesz-Thorin interpolation theorem, we finish the proof of {\bf Theorem One}.

\section{Point-wise estimate on partial operators}
\setcounter{equation}{0}
In order to prove {\bf Proposition One}, we consider a generalized fractional integral operator
\bel{I_alpha beta theta}
\I_{\alpha\beta\theta} f(x)~=~\int_{\R^3} f(y) \V^{\alpha\beta\vartheta}(x-y)dy
\eeq
where
\bel{V theta}
\V^{\alpha\beta\vartheta}(x)~=~|x_1|^{\alpha-1}|x_2|^{\alpha-1}|x_3|^{\beta-1} \Bigg[ {|x_1||x_2|\over |x_3|}+{|x_3|\over |x_1||x_2|}\Bigg]^{-\theta},\qquad \theta>0
\eeq
for $x_i\neq0, i=1,2,3$. In addition to $-1<\alpha<1$ and $0<\alpha+\beta<2$, we require $-\theta<\beta<\theta$ if $\theta\leq1$ and $-1<\beta<1$ if $\theta>1$.

Let ${\alpha+\beta\over 2}={1\over p}-{1\over q},1<p<q<\infty$ for $q\in\Z$ satisfying $(q-2)\left[{\alpha+\beta\over 2}\right]=(q-2)\left[{1\over p}-{1\over q}\right]\ge1$. Clearly,  $q \left[{1\over p}-{1\over q}\right]>1$ implies ${1\over p}-{2\over q}>0$.  There exists a $\vartheta=\vartheta(\alpha,\beta,q)>0$ such that 
\bel{theta constraints}
\alpha-\vartheta~>~0,\qquad \beta-{1\over p}+\vartheta~=~{1\over p}-{2\over q}-\alpha+\vartheta~>~0\qquad\hbox{if}\qquad \alpha>0.
\eeq
Observe that the second inequality in (\ref{theta constraints}) holds for every $\vartheta>0$ if $\alpha\leq0$.

From now on, we set $\theta\ge\vartheta$ if $\alpha>0$ and $\theta>0$ if $\alpha\leq0$.
\begin{remark} Consequently, we have
$\beta-{1\over p}+\theta={1\over p}-{2\over q}-\alpha+\theta>0$.
\end{remark}
Let $\Delta_\ell \I_{\alpha\beta\theta}$, $\Delta_{\ell j} \I_{\alpha\beta\theta}$ and $\Delta_{\ell jk} \I_{\alpha\beta\theta}$ defined in analogue to (\ref{Delta_l I}) for every $\ell,j,l\in\Z$.

Recall $\Gamma_{\ell jk}(x)$ given in (\ref{Gamma_ljk}). From (\ref{V theta}), we find
\bel{V bounds}
\begin{array}{lr}\ds
\V^{\alpha\beta\vartheta}(x-y)~=~|x_1-y_1|^{\alpha-1}|x_2-y_2|^{\alpha-1}|x_3-y_3|^{\beta-1} \Bigg[ {|x_1-y_1||x_2-y_2|\over |x_3-y_3|}+{|x_3-y_3|\over |x_1-y_1||x_1-y_2|}\Bigg]^{-\theta}
\\\\ \ds~~~~~~~~~~~~~~~~~~~~
~\approx~~2^{j(\alpha-1)} 2^{(j-\ell)(\alpha-1)}2^{\big[j+(j-\ell)-k\big](\beta-1)} \Big[2^k+2^{-k}\Big]^{-\theta}
\end{array}
\eeq
whenever $y\in\Gamma_{\ell jk}(x)$. 
By using (\ref{V bounds}), we have
\bel{Est.1}
\begin{array}{lr}\ds
\int_{\Gamma_{\ell jk}(x)} f(y) \V^{\alpha\beta\theta}(x-y)dy
\\\\ \ds
~\lesssim~2^{j(\alpha-1)} 2^{(j-\ell)(\alpha-1)}2^{\big[j+(j-\ell)-k\big](\beta-1)} \Big[2^k+2^{-k}\Big]^{-\theta} \int_{\Gamma_{\ell jk}(x)} f(y)dy
\\\\ \ds
~=~2^{j\alpha} 2^{(j-\ell)\alpha}2^{\big[j+(j-\ell)-k\big]\beta} \Big[2^k+2^{-k}\Big]^{-\theta} \left\{{1\over 2^{j} 2^{(j-\ell)}2^{\big[j+(j-\ell)-k\big]} }\int_{\Gamma_{\ell jk}(x)} f(y)dy\right\}
\\\\ \ds
~\leq~2^{j\alpha} 2^{(j-\ell)\alpha}2^{\big[j+(j-\ell)-k\big]\beta} \Big[2^k+2^{-k}\Big]^{-\theta}  \M f(x)
\end{array}
\eeq
where $\M f$ is the strong maximal function defined in $\R^3$.

From (\ref{Est.1}), we find
\bel{Est.2}
\begin{array}{lr}\ds
\int_{\Gamma_{\ell j}(x)} f(y) \V^{\alpha\beta\theta}(x-y)dy
~=~\sum_{k\in\Z} \int_{\Gamma_{\ell jk}(x)} f(y) \V^{\alpha\beta\theta}(x-y)dy
\\\\ \ds~~~~~~~~~~~~~~~~~
~\lesssim~2^{j\alpha} 2^{(j-\ell)\alpha}2^{\big[j+(j-\ell)\big]\beta}\left\{\sum_{k\in\Z}2^{-k\beta} \Big[2^k+2^{-k}\Big]^{-\theta}\right\}  \M f(x)
\\\\ \ds~~~~~~~~~~~~~~~~~
~\leq~2^{j\alpha} 2^{(j-\ell)\alpha}2^{\big[j+(j-\ell)\big]\beta}~ \M f(x)\left\{\begin{array}{lr}\ds\sum_{k>0} 2^{-k\theta}+\sum_{k\leq0} 2^{k(\theta-\beta)},\qquad \hbox{\small{$0<\beta<\theta$}}
\\\\ \ds
\sum_{k>0} 2^{-k(\beta+\theta)}+\sum_{k\leq0} 2^{k\theta},\qquad \hbox{\small{$-\theta<\beta\leq0$}}
\end{array}
\right.  
\\\\ \ds~~~~~~~~~~~~~~~~~
~\leq~\C_{\alpha~\beta~\theta}~ 2^{\big[j+(j-\ell)\big](\alpha+\beta)}~ \M f(x).
\end{array}
\eeq
On the other hand, by using H\"{o}lder inequality, we have
\bel{Est.3}
\begin{array}{lr}\ds
\int_{\Gamma_{\ell jk}(x)} f(y) \V^{\alpha\beta\theta}(x-y)dy
~\leq~\left\|f\right\|_{\L^p\left(\Gamma_{\ell jk}(x)\right)} \left\{ \int_{\Gamma_{\ell jk}(x)} \left[\V^{\alpha\beta\theta}(x-y)\right]^{p\over p-1} dy\right\}^{p-1\over p}
\\\\ \ds
~\lesssim~\left\|f\right\|_{\L^p\big(\Gamma_\ell(x)\big)}~2^{j(\alpha-1)} 2^{(j-\ell)(\alpha-1)}2^{\big[j+(j-\ell)-k\big](\beta-1)} \Big[2^k+2^{-k}\Big]^{-\theta} \left[2^j 2^{j-\ell} 2^{\left[j+(j-\ell)-k\right]}\right]^{1-{1\over p}} 
\\\\ \ds
~=~\left\|f\right\|_{\L^p\big(\Gamma_\ell(x)\big)}~2^{j\big(\alpha-{1\over p}\big)} 2^{(j-\ell)\big(\alpha-{1\over p}\big)}2^{\big[j+(j-\ell)-k\big]\big(\beta-{1\over p}\big)} \Big[2^k+2^{-k}\Big]^{-\theta}.
\end{array}
\eeq
Recall {\bf Remark 3.1}. From (\ref{Est.3}), we find
\bel{Est.4}
\begin{array}{lr}\ds
\int_{\Gamma_{\ell j}(x)} f(y) \V^{\alpha\beta\theta}(x-y)dy
~=~\sum_{k\in\Z} \int_{\Gamma_{\ell j k}(x)} f(y) \V^{\alpha\beta\theta}(x-y)dy
\\\\ \ds
~\lesssim~\left\|f\right\|_{\L^p\big(\Gamma_\ell(x) \big)}~2^{j\big(\alpha-{1\over p}\big)} 2^{(j-\ell)\big(\alpha-{1\over p}\big)}2^{\big[j+(j-\ell)\big]\big(\beta-{1\over p}\big)} \sum_{k\in\Z} 2^{-k\big(\beta-{1\over p}\big)}\Big[2^k+2^{-k}\Big]^{-\theta}
\\\\ \ds
~\leq~\left\|f\right\|_{\L^p\big(\Gamma_\ell(x) \big)}~2^{\big[j+(j-\ell)\big]\big(\alpha+\beta-{2\over p}\big)}
\left\{\begin{array}{lr}\ds
\sum_{k\leq0} 2^{k\theta}+\sum_{k>0} 2^{-k\big[\beta-{1\over p}+\theta\big]},\qquad ~~\hbox{\small{$\beta-{1\over p}\leq0$}}
\\\\ \ds
\sum_{k\leq0} 2^{-k\big[\beta-{1\over p}-\theta\big]}+\sum_{k>0} 2^{-k\theta},\qquad \hbox{\small{$\beta-{1\over p}>0$}}
\end{array}\right.
\\\\ \ds
~\leq~\C_{\alpha~\beta~\theta}~\left\|f\right\|_{\L^p\big(\Gamma_\ell(x)\big)}~2^{\big[j+(j-\ell)\big]\big(\alpha+\beta-{2\over p}\big)}.
\end{array}
\eeq
Given a non-zero function $f\in\L^p(\R^3)$, we define 
\bel{theta_l}
\varphi_\ell(x)~=~\left\| f\right\|_{\L^p(\R^3)}^{-p}\int_{\Gamma_\ell(x)} \Big(f(y)\Big)^p dy,\qquad \ell\in\Z. 
\eeq
\begin{remark} Clearly,  $0\leq\varphi_\ell(x)<1$ and $\sum_{\ell\in\Z} \varphi_\ell(x)=1$ for every $x\in\R^3$. 
\end{remark}
Next, we define $\lambda(\ell,x)\in\R$ implicitly by requiring
\bel{lambda}
\left[ {\varphi_\ell(x)\over \Big(\M f\Big)^p(x)}\left\| f\right\|_{\L^p(\R^3)}^p\right]^{1\over 2}~=~2^{\lambda(\ell,x)}2^{\lambda(\ell,x)-\ell}.
\eeq
For $-\infty<j\leq\lambda(\ell,x)$, by inserting (\ref{lambda}) to (\ref{Est.2}), we find
\bel{Est.5}
\begin{array}{lr}\ds
\Delta_{\ell j} \I_{\alpha\beta\theta} f(x)~\leq~\C_{\alpha~\beta~\theta}~ 2^{\big[j+(j-\ell)\big](\alpha+\beta)}~ \M f(x)
\\\\ \ds~~~~~~~~~~~~~~~~~~
~=~\C_{\alpha~\beta~\theta}~2^{2\big[j-\lambda(\ell,x)\big](\alpha+\beta)} ~2^{\big[\lambda(\ell,x)+\lambda(\ell,x)-\ell\big](\alpha+\beta)}~ \M f(x)
\\\\ \ds~~~~~~~~~~~~~~~~~~
~=~\C_{\alpha~\beta~\theta}~2^{2\big[j-\lambda(\ell,x)\big](\alpha+\beta)} ~\left[ {\varphi_\ell(x)\over \Big(\M f\Big)^p(x)}\left\| f\right\|_{\L^p(\R^3)}^p\right]^{\alpha+\beta\over 2}~\M f(x)
\\\\ \ds~~~~~~~~~~~~~~~~~~
~=~\C_{\alpha~\beta~\theta}~2^{2\big[j-\lambda(\ell,x)\big](\alpha+\beta)} ~\left[ {\varphi_\ell(x)\over \Big(\M f\Big)^p(x)}\left\| f\right\|_{\L^p(\R^3)}^p\right]^{{1\over p}-{1\over q}}~\M f(x)
\\\\ \ds~~~~~~~~~~~~~~~~~~
~=~\C_{\alpha~\beta~\theta}~2^{2\big[j-\lambda(\ell,x)\big](\alpha+\beta)}~\left\| f\right\|_{\L^p(\R^3)}^{1-{p\over q}} \Big(\vartheta_\ell(x)\Big)^{{1\over p}-{1\over q}} ~\Big(\M f\Big)^{p\over q}(x).
\end{array}
\eeq
For $\lambda(\ell,x)<j<\infty$, by inserting (\ref{lambda}) to (\ref{Est.4}), we have
\bel{Est.6}
\begin{array}{lr}\ds
\Delta_{\ell j} \I_{\alpha\beta\theta} f(x)~\leq~\C_{\alpha~\beta~\theta}~ \left\|f\right\|_{\L^p\big(\Gamma_\ell(x)\big)}~2^{\big[j+(j-\ell)\big]\big(\alpha+\beta-{2\over p}\big)}
\\\\ \ds~~~~~~~~~~~~~~~~~
~=~\C_{\alpha~\beta~\theta}~ \left\|f\right\|_{\L^p(\R^3)}~\Big(\varphi_\ell(x)\Big)^{1\over p}~2^{\big[j+(j-\ell)\big]\big(\alpha+\beta-{2\over p}\big)}
\\\\ \ds~~~~~~~~~~~~~~~~~
~=~\C_{\alpha~\beta~\theta}~ \left\|f\right\|_{\L^p(\R^3)}~\Big(\varphi_\ell(x)\Big)^{1\over p}~2^{\big[\lambda(\ell,x)+\lambda(\ell,x)-\ell\big]\big(\alpha+\beta-{2\over p}\big)}~2^{2\big[j-\lambda(\ell,x)\big]\big(\alpha+\beta-{2\over p}\big)}
\\\\ \ds~~~~~~~~~~~~~~~~~
~=~\C_{\alpha~\beta~\theta}~ \left\|f\right\|_{\L^p(\R^3)}~\Big(\varphi_\ell(x)\Big)^{1\over p}~
\left[ {\varphi_\ell(x)\over \Big(\M f\Big)^p(x)}\left\| f\right\|_{\L^p(\R^3)}^p\right]^{{\alpha+\beta\over 2}-{1\over p}}~2^{2\big[j-\lambda(\ell,x)\big]\big(\alpha+\beta-{2\over p}\big)}
\\\\ \ds~~~~~~~~~~~~~~~~~
~=~\C_{\alpha~\beta~\theta}~ \left\|f\right\|_{\L^p(\R^3)}~\Big(\varphi_\ell(x)\Big)^{1\over p}~
\left[ {\varphi_\ell(x)\over \Big(\M f\Big)^p(x)}\left\| f\right\|_{\L^p(\R^3)}^p\right]^{-{1\over q}}~2^{2\big[j-\lambda(\ell,x)\big]\big(\alpha+\beta-{2\over p}\big)}
\\\\ \ds~~~~~~~~~~~~~~~~~
~=~\C_{\alpha~\beta~\theta}~2^{2\big[j-\lambda(\ell,x)\big]\big(\alpha+\beta-{2\over p}\big)}~ \left\|f\right\|_{\L^p(\R^3)}^{1-{p\over q}}~\Big(\varphi_\ell(x)\Big)^{{1\over p}-{1\over q}} \Big(\M f\Big)^{p\over q}(x)
\\\\ \ds~~~~~~~~~~~~~~~~~
~=~\C_{\alpha~\beta~\theta}~2^{2\big[j-\lambda(\ell,x)\big]\big(-{2\over q}\big)}~ \left\|f\right\|_{\L^p(\R^3)}^{1-{p\over q}}~\Big(\varphi_\ell(x)\Big)^{{1\over p}-{1\over q}} \Big(\M f\Big)^{p\over q}(x).
\end{array}
\eeq
Denote
\bel{sigma}
\sigma~=~\min\left\{\alpha+\beta, {2\over q}\right\}.
\eeq
By putting together (\ref{Est.5}) and (\ref{Est.6}), we find
\bel{Est.7}
\Delta_{\ell j}  \I_{\alpha\beta\theta} f(x)~\leq~\C_{\alpha~\beta~\theta}~2^{-2\left| j-\lambda(\ell,x)\right| \sigma} \left\|f\right\|_{\L^p(\R^3)}^{1-{p\over q}}~\Big(\varphi_\ell(x)\Big)^{{1\over p}-{1\over q}} \Big(\M f\Big)^{p\over q}(x).
\eeq

By using (\ref{Est.7}) and summing over every $j\in\Z$, we obtain
\bel{Est.8}
\Delta_\ell  \I_{\alpha\beta\theta} f(x)~\leq~\C_{\alpha~\beta~\theta~q}~ \left\|f\right\|_{\L^p(\R^3)}^{1-{p\over q}}~\Big(\varphi_\ell(x)\Big)^{{1\over p}-{1\over q}} \Big(\M f\Big)^{p\over q}(x).
\eeq
Recall $(q-2)\left[{1\over p}-{1\over q}\right]\ge1$. From (\ref{Est.8}), we have
\bel{Est.9}
\begin{array}{lr}\ds
\int_{\R^3} \sum_{\ell\in\Z} \Delta_\ell \I_{\alpha\beta\theta} f(x) \Big(\Delta_{\ell-h}\I_{\alpha\beta} f\Big)^{q-1}(x)dx
\\\\ \ds
~\leq~\C_{\alpha~\beta~\theta~q}~ \left\|f\right\|_{\L^p(\R^3)}^{q-p}~\int_{\R^3}\Big(\M f\Big)^p(x)\sum_{\ell\in\Z}\Big(\varphi_\ell(x)\Big)^{{1\over p}-{1\over q}} \Big(\varphi_{\ell-h}(x)\Big)^{(q-1)\big[{1\over p}-{1\over q}\big]} dx
\\\\ \ds
~\leq~\C_{\alpha~\beta~\theta~q}~ \left\|f\right\|_{\L^p(\R^3)}^{q-p}~\int_{\R^3}\Big(\M f\Big)^p(x) dx\qquad \hbox{\small{by {\bf Remark 3.2}}}
\\\\ \ds
~\leq~\C_{\alpha~\beta~\theta~q}~ \left\|f\right\|_{\L^p(\R^3)}^q.
\end{array}
\eeq
\v

{\bf Proposition Two}:~~{\it Let ${\alpha+\beta\over 2}={1\over p}-{1\over q},1<p<q<\infty$. Suppose $q\in\Z$ satisfying $(q-2)\left[{\alpha+\beta\over 2}\right]\ge1$. Moreover, $\alpha>0$ and  $\vartheta=\vartheta(\alpha,\beta)>0$ is implicitly defined in (\ref{theta constraints}). We have
\bel{Lemma Two}
\int_{\R^3} \sum_{\ell\in\Z} \Delta_\ell \I_{\alpha\beta\vartheta} f(x) \Big(\Delta_{\ell-h}\I_{\alpha\beta} f\Big)^{q-1}(x)dx~\leq~\C_{\alpha~\beta~q} ~2^{-\ve |h|}~\left\| f\right\|_{\L^p(\R^3)}^q
\eeq
for some $\ve=\ve(\alpha,\beta,q)>0$.}

\section{Some implications on Proposition One}
\setcounter{equation}{0}

\subsection{Proposition Two implies Proposition One}
For $\alpha>0$, there exists a $\vartheta>0$ satisfying (\ref{theta constraints}). Let $\V^{\alpha\beta\theta}(x)$ defined in (\ref{V theta}). We find 
\bel{V compara}
\begin{array}{lr}\ds
\V^{\alpha\beta}(x)~=~|x_1|^{\alpha-1}|x_2|^{\alpha-1}|x_3|^{\beta-1} \Bigg[ {|x_1||x_2|\over |x_3|}+{|x_3|\over |x_1||x_2|}\Bigg]^{-1}
\\\\ \ds~~~~~~~~~~~
~\leq~|x_1|^{\alpha-1}|x_2|^{\alpha-1}|x_3|^{\beta-1} \Bigg[ {|x_1||x_2|\over |x_3|}+{|x_3|\over |x_1||x_2|}\Bigg]^{-\vartheta}~=~\V^{\alpha\beta\vartheta(x)}
\end{array}
\eeq
for every $x_i\neq0,i=1,2,3$. This implies $\Delta_\ell \I_{\alpha\beta} f(x)\leq \Delta_\ell \I_{\alpha\beta\vartheta} f(x)$. Therefore, {\bf Proposition Two} implies {\bf Proposition One}.

Suppose $\alpha\leq0$. We aim to find $-1<\alpha_1<\alpha\leq0<\alpha_2<1$ and $-1<\beta_1,\beta_2<1$ such that ${\alpha+\beta\over 2}={1\over p}-{1\over q}={\alpha_1+\beta_1\over 2}={\alpha_2+\beta_2\over 2}$. 
Choose $\alpha_1=\alpha/2$ and $\alpha_2={1\over 2}$. Simultaneously, $\beta_1, \beta_2$ will be fixed. There exists a $0<t<1$ such that $\alpha=(1-t)\alpha_1+t\alpha_2$ and  $\beta=(1-t)\beta_1+t\beta_2$.
Because $\alpha_2>\beta_2$,  there is a $\vartheta=\vartheta(\alpha_2,\beta_2,q)>0$ satisfying (\ref{theta constraints})  and $-\vartheta<\beta_2<\vartheta$
such that 
\bel{formula t}
 1~=~(1-t)\theta+t\vartheta\qquad\hbox{for some}\qquad \theta>1.
\eeq
\begin{remark} From the above, the value of $\alpha_i,\beta_i, i=1,2$ and $\theta$ depend on $\alpha,\beta,q$.
\end{remark}
Let $z\in\Cx$ and $0\leq \Re z\leq1$. Consider
\bel{I_z}
\begin{array}{lr}\ds
\Delta_\ell \I_{(1-z)\alpha_1+z\alpha_2~(1-z)\beta_1+z\beta_2~ (1-z)\theta+z\vartheta} f(x)~=~
\\\\ \ds
\int_{\Gamma_\ell(x)}f(y) \V^{(1-z)\alpha_1+z\alpha_2~(1-z)\beta_1+z\beta_2~ (1-z)\theta+z\vartheta}(x-y)dy
\end{array}
\eeq
where
\bel{V z}
\begin{array}{lr}\ds
\V^{(1-z)\alpha_1+z\alpha_2~(1-z)\beta_1+z\beta_2~ (1-z)\theta+z\vartheta}(x)~=~
\\\\ \ds
|x_1|^{(1-z)\alpha_1+z\alpha_2-1}|x_2|^{(1-z)\alpha_1+z\alpha_2-1}|x_3|^{(1-z)\beta_1+z\beta_2-1} \Bigg[ {|x_1||x_2|\over |x_3|}+{|x_3|\over |x_1||x_2|}\Bigg]^{-\big[(1-z)\theta+z\vartheta\big]}
\end{array}
\eeq
for $x_i\neq0, i=1,2,3$.

We have
\bel{V norm compara}
\left| \V^{(1-z)\alpha_1+z\alpha_2~(1-z)\beta_1+z\beta_2~ (1-z)\theta+z\vartheta}(x)\right|~=~\V^{(1-\Re z)\alpha_1+\Re z\alpha_2~(1-\Re z)\beta_1+\Re z\beta_2~ (1-\Re z)\theta+\Re z\vartheta}(x).
\eeq
Define
\bel{U_z}
\U(z)~=~\int_{\R^3} \sum_{\ell\in\Z} \Delta_\ell \I_{(1-z)\alpha_1+z\alpha_2~(1-z)\beta_1+z\beta_2~ (1-z)\theta+z\vartheta} f(x) \Big(\Delta_{\ell-h}\I_{\alpha\beta} f\Big)^{q-1}(x)dx.
\eeq
Observe that $\U(z)$ is analytic in the strip: $0\leq\Re z\leq1$ provided that
\bel{U z norm}
\begin{array}{lr}\ds
\left|\U(z)\right|~\leq~\int_{\R^3} \sum_{\ell\in\Z} \Delta_\ell \I_{(1-\Re z)\alpha_1+\Re z\alpha_2~(1-\Re z)\beta_1+\Re z\beta_2~ (1-\Re z)\theta+\Re z\vartheta} f(x) \Big(\Delta_{\ell-h}\I_{\alpha\beta} f\Big)^{q-1}(x)dx
\\\\ \ds~~~~~~~~~~
~\leq~\C_{\Re z~\alpha_1~\beta_1~\alpha_2~\beta_2~\theta~\vartheta~q}~\left\| f\right\|_{\L^p(\R^3)}^q\qquad\hbox{\small{by (\ref{Est.9})}}.
\end{array}
\eeq
Moreover, we have
\bel{U 0}
\left|\U(0+\i \Im z)\right|~\leq~\C_{\alpha_1~\beta_1~\theta~q}~\left\| f\right\|_{\L^p(\R^3)}^q.
\eeq
On the other hand, by applying {\bf Proposition Two}, we find
\bel{U 1}
\left|\U(1+\i \Im z)\right|~\leq~\C_{\alpha_2~\beta_2~\vartheta~q}~2^{-\ve|h|}\left\| f\right\|_{\L^p(\R^3)}^q.
\eeq
for some $\ve=\ve(\alpha_2,\beta_2,q)>0$.

Recall {\bf Remark 4.1}. From (\ref{U 0}) and (\ref{U 1}), by applying Three-Line lemma, we obtain
\bel{U t}
\begin{array}{lr}\ds
\U(t)~=~\int_{\R^3} \sum_{\ell\in\Z} \Delta_\ell \I_{\alpha\beta} f(x) \Big(\Delta_{\ell-h}\I_{\alpha\beta} f\Big)^{q-1}(x)dx
~\leq~\C_{\alpha~\beta~q} ~2^{-t\ve|h|}~\left\| f\right\|_{\L^p(\R^3)}^q.
\end{array}
\eeq

\subsection{Simplifying Proposition Two}
In order to prove {\bf Proposition Two}, we aim to show
\bel{almost Ortho Est.1}
\begin{array}{lr}\ds
\int_{\R^3}  \Delta_\ell \I_{\alpha\beta\vartheta} f(x) \Big(\Delta_{\ell-h}\I_{\alpha\beta} f\Big)^{q-1}(x)dx
\\\\ \ds
~\leq~\C_{\alpha~\beta~q}~2^{-\ve|h|}~\left\|f\right\|_{\L^p(\R^3)}^{q-p}~\int_{\R^3}\Big(\varphi_{\ell-h}(x)\Big)^{(q-2)\big[{1\over p}-{1\over q}\big]} \Big(\M f\Big)^p(x) dx
\end{array}
\eeq
for some $\ve=\ve(\alpha,\beta,q)>0$. 

Recall  $\varphi_\ell(x)$ defined in (\ref{theta_l}) and {\bf Remark 3.2}.
We have 
$\ds\sum_{\ell\in\Z}\Big(\varphi_\ell(x)\Big)^{(q-2)\big[{1\over p}-{1\over q}\big]}\leq1$ provided by $(q-2)\left[{1\over p}-{1\over q}\right]\ge1$. From (\ref{almost Ortho Est.1}), by summing all $\ell\in\Z$ and using the $\L^p$-boundedness of $\M$, 
we obtain (\ref{Lemma}) in {\bf Proposition One}.

Next, we claim that it is suffice to prove (\ref{almost Ortho Est.1}) for $\ell=h$. Namely, 
\bel{almost Ortho Est.2}
\begin{array}{lr}\ds
\int_{\R^3}  \Delta_h \I_{\alpha\beta\vartheta} f(x) \Big(\Delta_0\I_{\alpha\beta} f\Big)^{q-1}(x)dx
\\\\ \ds
~\leq~\C_{\alpha~\beta~q}~2^{-\ve|h|}~\left\|f\right\|_{\L^p(\R^3)}^{q-p}~\int_{\R^3}\Big(\varphi_0(x)\Big)^{(q-2)\big[{1\over p}-{1\over q}\big]} \Big(\M f\Big)^p(x) dx.
\end{array}
\eeq
Denote 
\bel{tau}
\tau_s x~=~\Big(x_1, 2^{-s}x_2, 2^{-s} x_3\Big),\qquad s\in\Z.
\eeq
Let $\V^{\alpha\beta\theta}(x)$ defined in (\ref{V theta}). We find
\bel{V dila}
\V^{\alpha\beta\vartheta}(\tau_s x)~=~2^{s(1-\alpha)}2^{s(1-\beta)} \V^{\alpha\beta\vartheta}(x)
\eeq
for $x_i\neq0, i=1,2,3$.

Recall $\Gamma_\ell(x)$ defined in (\ref{Gamma_ljk}). From (\ref{tau})-(\ref{V dila}), we find
\bel{almost Ortho dila}
\begin{array}{lr}\ds
\Delta_\ell \I_{\alpha\beta\vartheta} f(\tau_s x)~=~\int_{\Gamma_\ell(\tau_s x)} f(y)\V^{\alpha\beta\vartheta}(\tau_s x-y)dy
\\\\ \ds~~~~~~~~~~~~~~~~~~~~~
~=~\int_{\Gamma_{\ell-s}(x)} f(\tau_s y)\V^{\alpha\beta\vartheta}(\tau_s x-\tau_s y)2^{-s}2^{-s}dy
\\\\ \ds~~~~~~~~~~~~~~~~~~~~~
~=~2^{-s (\alpha+\beta)}\int_{\Gamma_{\ell-s}(x)} f(\tau_s y)\V^{\alpha\beta\vartheta}(x-y)dy.
\end{array}
\eeq
Let $f_s(x)=f(\tau_s x)$. Clearly, we have 
\bel{f norm s}
\left\| f_s\right\|_{\L^p(\R^3)}^p~=~2^s2^s \left\| f\right\|_{\L^p(\R^3)}^p.
\eeq
Note that $\varphi_\ell(x)$ defined in (\ref{theta_l}) depends on the given function
$f\in\L^p(\R^3)$. We write
\bel{theta dila}
\begin{array}{lr}\ds
\varphi_\ell^s (x)~\doteq~\left\| f_s\right\|_{\L^p(\R^3)}^{-p}\int_{\Gamma_\ell( x)} \Big(f_s(y)\Big)^p dy
\\\\ \ds~~~~~~~~
~=~\left\| f_s\right\|_{\L^p(\R^3)}^{-p}2^{s}2^{s}\int_{\Gamma_\ell( x)} \Big(f(\tau_s y)\Big)^p d(\tau_s y)
\\\\ \ds~~~~~~~~
~=~\left\| f\right\|_{\L^p(\R^3)}^{-p}\int_{\Gamma_{\ell+s}(\tau_s x)} \Big(f( y)\Big)^p dy
\\\\ \ds~~~~~~~~
~=~\varphi_{\ell+s}(\tau_s x).
\end{array}
\eeq
Furthermore, a direct computation shows
\bel{Mf s}
\M f_s(x)~=~\M f(\tau_s x).
\eeq
Suppose (\ref{almost Ortho Est.2}) hold. 
By using (\ref{almost Ortho dila})-(\ref{Mf s}) and taking into account for $s=\ell-h$, we have
\bel{almost Ortho transform}
\begin{array}{lr}\ds
\int_{\R^3}  \Delta_\ell \I_{\alpha\beta\vartheta} f(x) \Big(\Delta_{\ell-h}\I_{\alpha\beta} f\Big)^{q-1}(x)dx
\\\\ \ds
~=~2^{-(\ell-h)}2^{-(\ell-h)}\int_{\R^3}  \Delta_\ell \I_{\alpha\beta\vartheta} f(\tau_{\ell-h} x) \Big(\Delta_{\ell-h}\I_{\alpha\beta} f\Big)^{q-1}(\tau_{\ell-h} x)dx
\\\\ \ds
~=~2^{-(\ell-h)}2^{-(\ell-h)} 2^{-(\ell-h)(\alpha+\beta)q} \int_{\R^3}  \Delta_h \I_{\alpha\beta\vartheta} f_{\ell-h}(x) \Big(\Delta_0\I_{\alpha\beta} f_{\ell-h}\Big)^{q-1}(x)dx\qquad\hbox{\small{by (\ref{almost Ortho dila})}}
\\\\ \ds
~\leq~\C_{\alpha~\beta~q} 2^{-(\ell-h)}2^{-(\ell-h)} 2^{-(\ell-h)(\alpha+\beta)q}  
\\\\ \ds~~~~~~~
2^{-\ve|h|}~\left\|f_{\ell-h}\right\|_{\L^p(\R^3)}^{q-p}~\int_{\R^3}\Big(\varphi_0^{\ell-h}(x)\Big)^{(q-2)\big[{1\over p}-{1\over q}\big]} \Big(\M f_{\ell-h}\Big)^p(x) dx \qquad\hbox{\small{by (\ref{almost Ortho Est.2})}}
\\\\ \ds
~=~\C_{\alpha~\beta~q}~2^{-(\ell-h)(\alpha+\beta)q} 2^{(\ell-h)\big[{q\over p}-1\big]} 2^{(\ell-h)\big[{q\over p}-1\big]} 
\\\\ \ds~~~~~~~
2^{-\ve|h|}~\left\|f\right\|_{\L^p(\R^3)}^{q-p}~\int_{\R^3}\Big(\varphi_{\ell-h}(\tau_{\ell-h}x)\Big)^{(q-2)\big[{1\over p}-{1\over q}\big]} \Big(\M f\Big)^p(\tau_{\ell-h}x) d\tau_{\ell-h}x\qquad \hbox{\small{by (\ref{theta dila})-(\ref{Mf s})}}
\\\\ \ds
~=~\C_{\alpha~\beta~q}~
2^{-\ve|h|}~\left\|f\right\|_{\L^p(\R^3)}^{q-p}~\int_{\R^3}\Big(\varphi_{\ell-h}(x)\Big)^{(q-2)\big[{1\over p}-{1\over q}\big]} \Big(\M f\Big)^p(x) dx
\end{array}
\eeq
where  $(\alpha+\beta)q=2\big[{q\over p}-1\big]$.

\begin{remark}
Recall $\Gamma_\ell(x)$ defined in (\ref{Gamma_ljk}) for $\ell\in\Z$. Observe that the union $\Gamma_h(x)=\Cup_{j\in\Z}\Gamma_{h j}(x)$ takes over every $j\in\Z$. If $h<0$, it is equivalent to define $\Gamma_h(x)$ for $h>0$ by switching the roles of $x_1-y_1$ and $x_2-y_2$. Because $\V^{\alpha\beta\vartheta}(x)$ defined in (\ref{V theta}) is symmetric $w.r.t$ $x_1$ and $x_2$, every regarding estimate remains the same.
For this symmetry reason, we prove (\ref{almost Ortho Est.2}) for $h>0$ only.
\end{remark}

\section{Proof of Proposition Two}
\setcounter{equation}{0}
Let ${\alpha+\beta\over 2}={1\over p}-{1\over q},1<p<q<\infty$ and $q\in\Z$ satisfying $(q-2)\left[{1\over p}-{1\over q}\right]$. Suppose $\alpha>0$ and $\vartheta=\vartheta(\alpha,\beta,q)>0$ is implicitly defined in (\ref{theta constraints}). In particular, we have $\alpha-\vartheta>0$.

From the previous section, we left to show
\bel{Lemma>}\begin{array}{lr}\ds
\int_{\R^3}  \Delta_\ell \I_{\alpha\beta\vartheta} f(x) \Big(\Delta_0\I_{\alpha\beta} f\Big)^{q-1}(x)dx
\\\\ \ds
~\leq~\C_{\alpha~\beta~q}~2^{-\ve\ell}~\left\|f\right\|_{\L^p(\R^3)}^{q-p}~\int_{\R^3}\Big(\varphi_0(x)\Big)^{(q-2)\big[{1\over p}-{1\over q}\big]} \Big(\M f\Big)^p(x) dx,\qquad \ell>0
\end{array}
\eeq
for some $\ve=\ve(\alpha,\beta,q)>0$.

First,  we write
\bel{Ortho Est.1}
\begin{array}{lr}\ds
\int_{\R^3}  \Delta_\ell \I_{\alpha\beta\vartheta} f(x) \Big(\Delta_0\I_{\alpha\beta} f\Big)^{q-1}(x)dx
\\\\ \ds
~=~ \int_{\R^3} \left\{\int_{\Gamma_\ell(x)} f(y)\V^{\alpha\beta\vartheta}(x-y)dy\right\}\prod_{m=1}^{q-1} \left\{\int_{\Gamma_0(x)} f(y)\V^{\alpha\beta}(x-y)dy\right\} dx
\\\\ \ds
~=~ \int_{\R^3} \sum_{j,j_1,\ldots, j_{q-1}\in\Z}\left\{\int_{\Gamma_{\ell j}(x)} f(y)\V^{\alpha\beta\vartheta}(x-y)dy\right\}\prod_{m=1}^{q-1} \left\{\int_{\Gamma_{0 j_m}(x)} f(y)\V^{\alpha\beta}(x-y)dy\right\} dx.
\end{array}
\eeq
Let 
$j_\nu=\min\big\{j_m, m=1,2,\ldots,q-1\big\}$.
 We develop a $3$-fold estimate $w.r.t$
\bel{Sum split 3}
\begin{array}{cc}\ds
\sum_{j,j_1,\ldots,j_{q-1}\in\Z}~=~\sum_{\G_1}+\sum_{\G_2}+\sum_{\G_3};
\\\\ \ds
\G_1~=~\left\{ j,j_1,\ldots,j_{q-1}\in\Z\colon j-\ell\ge j_\nu-2\right\},
\qquad
 \G_2~=~\left\{ j,j_1,\ldots,j_{q-1}\in\Z\colon j\leq j_\nu\right\},
\\\\ \ds
\G_3~=~\left\{ j,j_1,\ldots,j_{q-1}\in\Z\colon j-\ell< j_\nu-2<j-2\right\}.
\end{array}
\eeq
Denote $\jmath$ and $\jmath_m, m=1,2,\ldots,q-1$ implicitly by
\bel{jmath}
j~=~\lambda(\ell,x)+\jmath,\qquad  j_m~=~\lambda(0,x)+\jmath_m,\qquad m=1,2,\ldots,q-1.
\eeq
Recall $\varphi_\ell(x)$ defined in (\ref{theta_l})  and $\lambda(\ell,x)$ defined in (\ref{lambda}). We find
\bel{theta rewrite}
 \varphi_\ell(x)~=~{\Big(\M f\Big)^p(x)\over \left\| f\right\|_{\L^p(\R^3)}^p}~2^{2\lambda(\ell,x)}2^{2\lambda(\ell,x)-2\ell},\qquad\ell\in\Z.
\eeq

\subsection{Case 1: $j-\ell\ge j_\nu-2$}
Suppose $\lambda(\ell,x)-\lambda(0,x)>(1-\delta)\ell$ for some $0<\delta<{1\over 2}$. We have
\bel{theta Case1 est.1}
\begin{array}{lr}\ds
{\varphi_0(x)\over \varphi_\ell(x)}~=~2^{2\big[\lambda(0,x)-\lambda(\ell,x)\big]} ~2^{2\big[\lambda(0,x)-\lambda(\ell,x)\big]} 2^{2\ell}
~<~2^{-2(1-\delta)\ell} 2^{-2(1-\delta)\ell} 2^{2\ell}
~=~2^{-2(1-2\delta)\ell}.
\end{array}
\eeq
Recall the point-wise estimate in (\ref{Est.7}).  Denote $\sigma=\min\left\{\alpha+\beta, {2\over q}\right\}$. We find
\bel{Case1.Est.1}
\begin{array}{lr}\ds
\int_{\Gamma_{\ell j}(x)} f(y)\V^{\alpha\beta\vartheta}(x-y)dy\prod_{m=1}^{q-1} \int_{\Gamma_{0 j_m}(x)} f(y)\V^{\alpha\beta}(x-y)dy 
\\\\ \ds
~\leq~
\C_{\alpha~\beta~q}~2^{-2\left| j-\lambda(\ell,x)\right| \sigma} 
\prod_{m=1}^{q-1} 2^{-2\left| j_m-\lambda(0,x)\right| \sigma} 
\left\|f\right\|_{\L^p(\R^3)}^{q-p}~\Big(\varphi_\ell(x)\Big)^{{1\over p}-{1\over q}} \Big(\varphi_0(x)\Big)^{(q-1)\big[{1\over p}-{1\over q}\big]} \Big(\M f\Big)^p(x)
\\\\ \ds
~\leq~\C_{\alpha~\beta~q}~2^{-2\left| j-\lambda(\ell,x)\right| \sigma} 
\prod_{m=1}^{q-1} 2^{-2\left| j_m-\lambda(0,x)\right|\sigma} 
\\\\ \ds~~~~~~~
\left\|f\right\|_{\L^p(\R^3)}^{q-p}~2^{-2(1-2\delta)\big[{1\over p}-{1\over q}\big]\ell}
\Big(\varphi_\ell(x)\Big)^{2\big[{1\over p}-{1\over q}\big]} 
\Big(\varphi_0(x)\Big)^{(q-2)\big[{1\over p}-{1\over q}\big]} \Big(\M f\Big)^p(x)\qquad\hbox{\small{by (\ref{theta Case1 est.1})}}
\\\\ \ds
~\leq~\C_{\alpha~\beta~q}~2^{-2\left| j-\lambda(\ell,x)\right| \sigma} 
\prod_{m=1}^{q-1} 2^{-2\left| j_m-\lambda(0,x)\right|\sigma} 
\\\\ \ds~~~~~~~
\left\|f\right\|_{\L^p(\R^3)}^{q-p}~2^{-(1-2\delta)\big[{2\over p}-{2\over q}\big]\ell}
\Big(\varphi_0(x)\Big)^{(q-2)\big[{1\over p}-{1\over q}\big]} \Big(\M f\Big)^p(x)\qquad\hbox{\small{($0<\varphi_\ell<1,~\ell\in\Z$)}}
\end{array}
\eeq
On the other hand, suppose $\lambda(\ell,x)-\lambda(0,x)\leq (1-\delta)\ell$. As shown in (\ref{jmath}), we write 
$ j-\ell=\lambda(\ell,x)+\jmath-\ell$ and $\lambda(0,x)+\jmath_\nu-2=j_\nu-2$. 
Consequently, $j-\ell\ge j_\nu-2$ implies
\bel{Case1 index.est.1}
\jmath-\jmath_\nu~\ge~\ell-\big[ \lambda(\ell,x)-\lambda(0,x)\big]-2~\ge~\delta\ell -2.
\eeq
By using (\ref{Est.7}), we find
\bel{Case1.Est.2}
\begin{array}{lr}\ds
\int_{\Gamma_{\ell j}(x)} f(y)\V^{\alpha\beta\vartheta}(x-y)dy\prod_{m=1}^{q-1} \int_{\Gamma_{0 j_m}(x)} f(y)\V^{\alpha\beta}(x-y)dy 
\\\\ \ds
~\leq~
\C_{\alpha~\beta~q}~2^{-2\left| j-\lambda(\ell,x)\right| \sigma} 
\prod_{m=1}^{q-1} 2^{-2\left| j_m-\lambda(0,x)\right| \sigma} 
\left\|f\right\|_{\L^p(\R^3)}^{q-p}~\Big(\varphi_\ell(x)\Big)^{{1\over p}-{1\over q}} \Big(\varphi_0(x)\Big)^{(q-1)\big[{1\over p}-{1\over q}\big]} \Big(\M f\Big)^p(x)
\\\\ \ds
~\leq~\C_{\alpha~\beta~q}~2^{-2\left| \jmath\right| \sigma} 
\prod_{m=1}^{q-1} 2^{-2\left|\jmath_m\right| \sigma} 
\left\|f\right\|_{\L^p(\R^3)}^{q-p}~
\Big(\varphi_0(x)\Big)^{(q-2)\big[{1\over p}-{1\over q}\big]} \Big(\M f\Big)^p(x)
\qquad
\hbox{\small{by (\ref{jmath})}}
\\\\ \ds
~\leq~\C_{\alpha~\beta~q}~2^{-\left|\jmath-\jmath_\nu\right|\sigma} ~     2^{-\left| \jmath\right| \sigma} 
\prod_{m=1}^{q-1} 2^{-\left|\jmath_m\right| \sigma} 
\left\|f\right\|_{\L^p(\R^3)}^{q-p}~
\Big(\varphi_0(x)\Big)^{(q-2)\big[{1\over p}-{1\over q}\big]} \Big(\M f\Big)^p(x)
\\ \ds~~~~~~~~~~~~~~~~~~~~~~~~~~~~~~~~~~~~~~~~~~~~~~~~~~~~~~~~~~~~~~~~~~~~~~~~~~~~~~~~~~~~
\hbox{\small{( $|\jmath-\jmath_\nu|\leq|\jmath|+|\jmath_\nu|$ )}}
\\\\ \ds
~\leq~\C_{\alpha~\beta~q} ~     2^{-\left| j-\lambda(\ell,x)\right| \sigma} 
\prod_{m=1}^{q-1} 2^{-\left|j_m-\lambda(0,x)\right| \sigma} ~
2^{-\delta\sigma\ell}\left\|f\right\|_{\L^p(\R^3)}^{q-p}
\Big(\varphi_0(x)\Big)^{(q-2)\big[{1\over p}-{1\over q}\big]} \Big(\M f\Big)^p(x)
\\ \ds~~~~~~~~~~~~~~~~~~~~~~~~~~~~~~~~~~~~~~~~~~~~~~~~~~~~~~~~~~~~~~~~~~~~~~~~~~~~~~~~~~~~~~~~~~~~~~~~~~~~~~~~~~~~~~~~
\hbox{\small{by (\ref{Case1 index.est.1})}}.
\end{array}
\eeq
By putting together (\ref{Case1.Est.1}) and (\ref{Case1.Est.2}) with $\delta={1\over 3}$, we obtain
\bel{Case1.Est}
\begin{array}{lr}\ds
\sum_{\G_1}\int_{\Gamma_{\ell j}(x)} f(y)\V^{\alpha\beta\vartheta}(x-y)dy\prod_{m=1}^{q-1} \int_{\Gamma_{0 j_m}(x)} f(y)\V^{\alpha\beta}(x-y)dy
\\\\ \ds
~\leq~\C_{\alpha~\beta~q}~2^{-{1\over 3}\sigma\ell} \left\|f\right\|_{\L^p(\R^3)}^{q-p}
\Big(\varphi_0(x)\Big)^{(q-2)\big[{1\over p}-{1\over q}\big]} \Big(\M f\Big)^p(x)
\sum_{j,j_1,\ldots,j_{q-1}\in\Z}     2^{-\left| j-\lambda(\ell,x)\right| \sigma} 
\prod_{m=1}^{q-1} 2^{-\left|j_m-\lambda(0,x)\right|\sigma} 
\\\\ \ds
~\leq~\C_{\alpha~\beta~q}~2^{-{1\over 3}\sigma\ell} \left\|f\right\|_{\L^p(\R^3)}^{q-p}
\Big(\varphi_0(x)\Big)^{(q-2)\big[{1\over p}-{1\over q}\big]} \Big(\M f\Big)^p(x)
\end{array}
\eeq
where $\G_1=\left\{ j,j_1,\ldots,j_{q-1}\in\Z\colon j-\ell\ge j_\nu-2\right\}$.

\subsection{Case 2: $j\leq j_\nu$}
Suppose $\lambda(\ell,x)-\lambda(0,x)\leq\delta \ell$ for some $0<\delta<{1\over 2}$. We have
\bel{theta Case2 est.1}
\begin{array}{lr}\ds
{\varphi_\ell(x)\over \varphi_0(x)}~=~2^{2\big[\lambda(\ell,x)-\lambda(0,x)\big]} ~2^{2\big[\lambda(\ell,x)-\lambda(0,x)\big]} 2^{-2\ell}
~\leq~2^{2\delta\ell}2^{2\delta\ell}2^{-2\ell}~=~2^{-2(1-2\delta)\ell}.
\end{array}
\eeq
By using (\ref{Est.7}), we find
\bel{Case2.Est.1}
\begin{array}{lr}\ds
\int_{\Gamma_{\ell j}(x)} f(y)\V^{\alpha\beta\vartheta}(x-y)dy\prod_{m=1}^{q-1} \int_{\Gamma_{0 j_m}(x)} f(y)\V^{\alpha\beta}(x-y)dy 
\\\\ \ds
~\leq~
\C_{\alpha~\beta~q}~2^{-2\left| j-\lambda(\ell,x)\right| \sigma} 
\prod_{m=1}^{q-1} 2^{-2\left| j_m-\lambda(0,x)\right| \sigma} 
\left\|f\right\|_{\L^p(\R^3)}^{q-p}~\Big(\varphi_\ell(x)\Big)^{{1\over p}-{1\over q}} \Big(\varphi_0(x)\Big)^{(q-1)\big[{1\over p}-{1\over q}\big]} \Big(\M f\Big)^p(x)
\\\\ \ds
~\leq~\C_{\alpha~\beta~q}~2^{-2\left| j-\lambda(\ell,x)\right| \sigma} 
\prod_{m=1}^{q-1} 2^{-2\left| j_m-\lambda(0,x)\right|\sigma} 
\\ \ds~~~~~~~
\left\|f\right\|_{\L^p(\R^3)}^{q-p}~2^{-2(1-2\delta)\big[{1\over p}-{1\over q}\big]\ell}

\Big(\varphi_0(x)\Big)^{q\big[{1\over p}-{1\over q}\big]} \Big(\M f\Big)^p(x)\qquad\hbox{\small{by (\ref{theta Case2 est.1})}}
\\\\ \ds
~\leq~\C_{\alpha~\beta~q}~2^{-2\left| j-\lambda(\ell,x)\right| \sigma} 
\prod_{m=1}^{q-1} 2^{-2\left| j_m-\lambda(0,x)\right|\sigma} 
\\ \ds~~~~~~~
\left\|f\right\|_{\L^p(\R^3)}^{q-p}~2^{-(1-2\delta)\big[{2\over p}-{2\over q}\big]\ell}
\Big(\varphi_0(x)\Big)^{(q-2)\big[{1\over p}-{1\over q}\big]} \Big(\M f\Big)^p(x)\qquad\hbox{\small{by {\bf Remark 3.2}}}.
\end{array}
\eeq
On the other hand, suppose $\lambda(\ell,x)-\lambda(0,x)>\delta\ell$. As shown in (\ref{jmath}), we write 
$ j=\lambda(\ell,x)+\jmath$ and $\lambda(0,x)+\jmath_\nu=j_\nu$. 
Therefore, $j\leq j_\nu$ implies
\bel{Case2 index.est.1}
\jmath_\nu-\jmath~\ge~ \lambda(\ell,x)-\lambda(0,x)~>~\delta\ell.
\eeq
By using (\ref{Est.7}), we find
\bel{Case2.Est.2}
\begin{array}{lr}\ds
\int_{\Gamma_{\ell j}(x)} f(y)\V^{\alpha\beta\vartheta}(x-y)dy\prod_{m=1}^{q-1} \int_{\Gamma_{0 j_m}(x)} f(y)\V^{\alpha\beta}(x-y)dy 
\\\\ \ds
~\leq~
\C_{\alpha~\beta~}~2^{-2\left| j-\lambda(\ell,x)\right| \sigma} 
\prod_{m=1}^{q-1} 2^{-2\left| j_m-\lambda(0,x)\right| \sigma} 
\left\|f\right\|_{\L^p(\R^3)}^{q-p}~\Big(\varphi_\ell(x)\Big)^{{1\over p}-{1\over q}} \Big(\varphi_0(x)\Big)^{(q-1)\big[{1\over p}-{1\over q}\big]} \Big(\M f\Big)^p(x)
\\\\ \ds
~\leq~\C_{\alpha~\beta~q}~2^{-2\left| \jmath\right| \sigma} 
\prod_{m=1}^{q-1} 2^{-2\left|\jmath_m\right| \sigma} 
\left\|f\right\|_{\L^p(\R^3)}^{q-p}~
\Big(\varphi_0(x)\Big)^{(q-2)\big[{1\over p}-{1\over q}\big]} \Big(\M f\Big)^p(x)
\qquad
\hbox{\small{by (\ref{jmath})}}
\\\\ \ds
~\leq~\C_{\alpha~\beta~}~2^{-\left|\jmath-\jmath_\nu\right|\sigma} ~     2^{-\left| \jmath\right| \sigma} 
\prod_{m=1}^{q-1} 2^{-\left|\jmath_m\right| \sigma} 
\left\|f\right\|_{\L^p(\R^3)}^{q-p}~
\Big(\varphi_0(x)\Big)^{(q-2)\big[{1\over p}-{1\over q}\big]} \Big(\M f\Big)^p(x)
\\ \ds~~~~~~~~~~~~~~~~~~~~~~~~~~~~~~~~~~~~~~~~~~~~~~~~~~~~~~~~~~~~~~~~~~~~~~~~~~~~~~~~~~~~
\hbox{\small{( $|\jmath-\jmath_\nu|\leq|\jmath|+|\jmath_\nu|$ )}}
\\\\ \ds
~\leq~\C_{\alpha~\beta~} ~     2^{-\left| j-\lambda(\ell,x)\right| \sigma} 
\prod_{m=1}^{q-1} 2^{-\left|j_m-\lambda(0,x)\right| \sigma} ~
2^{-\delta\sigma\ell}\left\|f\right\|_{\L^p(\R^3)}^{q-p}
\Big(\varphi_0(x)\Big)^{(q-2)\big[{1\over p}-{1\over q}\big]} \Big(\M f\Big)^p(x)
\\ \ds~~~~~~~~~~~~~~~~~~~~~~~~~~~~~~~~~~~~~~~~~~~~~~~~~~~~~~~~~~~~~~~~~~~~~~~~~~~~~~~~~~~~~~~~~~~~~~~~~~~~~~~~~~~~~~~~
\hbox{\small{by (\ref{Case2 index.est.1})}}.
\end{array}
\eeq
By putting together  (\ref{Case2.Est.1}) and (\ref{Case2.Est.2}) with $\delta={1\over 3}$, we obtain
\bel{Case2.Est}
\begin{array}{lr}\ds
\sum_{\G_2}\int_{\Gamma_{\ell j}(x)} f(y)\V^{\alpha\beta\vartheta}(x-y)dy\prod_{m=1}^{q-1} \int_{\Gamma_{0 j_m}(x)} f(y)\V^{\alpha\beta}(x-y)dy
\\\\ \ds
~\leq~\C_{\alpha~\beta~q}~2^{-{1\over 3}\sigma\ell} \left\|f\right\|_{\L^p(\R^3)}^{q-p}
\Big(\varphi_0(x)\Big)^{(q-2)\big[{1\over p}-{1\over q}\big]} \Big(\M f\Big)^p(x)
\sum_{j,j_1,\ldots,j_{q-1}\in\Z}     2^{-\left| j-\lambda(\ell,x)\right| \sigma} 
\prod_{m=1}^{q-1} 2^{-\left|j_m-\lambda(0,x)\right|\sigma} 
\\\\ \ds
~\leq~\C_{\alpha~\beta~q}~2^{-{1\over 3}\sigma\ell} \left\|f\right\|_{\L^p(\R^3)}^{q-p}
\Big(\varphi_0(x)\Big)^{(q-2)\big[{1\over p}-{1\over q}\big]} \Big(\M f\Big)^p(x)
\end{array}
\eeq
where $\G_2=\left\{ j,j_1,\ldots,j_{q-1}\in\Z\colon j\leq j_\nu\right\}$.

\subsection{Case 3: $j-\ell < j_\nu-2<j-2$}
Recall $\Gamma_{\ell jk}(x)$ and $\Gamma_{\ell}(x), \Gamma_{\ell j}(x)$ defined in (\ref{Gamma_ljk}) and (\ref{Gamma_lj}) respectively.  We define \bel{Gamma*_ljk}
{^*}\Gamma_{\ell j k}(x)~=~\left\{ y\in\R^3\colon \left.\begin{array}{cc}\ds2^{j-3}\leq|x_1-y_1|<2^{j+3},~2^{j-3-\ell}\leq|x_2-y_2|<2^{j+3-\ell},
\\\\ \ds
2^{j+(j-\ell)-k}\leq|x_3-y_3|<2^{j+(j-\ell)+1-k}
\end{array}\right.\right\}.
\eeq
and
\bel{Gamma*_lj}
{^*}\Gamma_\ell(x)~=~\Cup_{j\in\Z}~ {^*}\Gamma_{\ell j}(x),\qquad {^*}\Gamma_{\ell j}(x)~=~\Cup_{k\in\Z}~ {^*}\Gamma_{\ell jk}(x).
\eeq
Furthermore, the projection of ${^*}\Gamma_{\ell j}(x)={^*}\Gamma_{\ell j}(x_1,x_2)$ in the $(x_1,x_2)$-subspace is denoted by
\bel{Gamma_ljk*}
\begin{array}{cc}\ds
{^*}\Lambda_{\ell j}(x_1,x_2)~=~{^*}\Lambda^1_{\ell j}(x_1){^*}\Lambda^2_{\ell j}(x_2),
\\\\ \ds
{^*}\Lambda^1_{\ell j}(x_1)~=~\left\{ y_1\in\R\colon 2^{j-3}\leq|x_1-y_1|<2^{j+3}\right\},
~~~~
 {^*}\Lambda^2_{\ell j}(x_2)~=~\left\{ y_2\in\R\colon 2^{j-3-\ell}\leq|x_2-y_2|<2^{j+3-\ell}\right\}.
\end{array}
\eeq
From direct computation, we have
\bel{Int change}
\begin{array}{lr}\ds
\int_{\R^3} \left\{\int_{\Gamma_{\ell j}(x)} f(y)\V^{\alpha\beta\vartheta}(x-y)dy\right\}\prod_{m=1}^{q-1} \left\{\int_{\Gamma_{0 j_m}(x)} f(y^m)\V^{\alpha\beta}(x-y^m)dy^m\right\} dx
\\\\ \ds
~=~\idotsint_{\R^3\times\cdots\times\R^3}f(y)\prod_{m=1}^{q-1} f(y^m)
\\ \ds~~~~~~
 \left\{\int_{\Gamma_{\ell j}(y)\cap \left[\Cap_{m=1}^{q-1} \Gamma_{0 j_m}(y^m)\right]}    \V^{\alpha\beta\vartheta}(x-y)\prod_{m=1}^{q-1} \V^{\alpha\beta}(x-y^m)dx\right\} dy \prod_{m=1}^{q-1} dy^m.
\end{array}
\eeq
\v
\begin{lemma} Let $r=j-j_\nu+2$. Suppose $\Gamma_{\ell j}(y)\cap \left[\Cap_{m=1}^{q-1} \Gamma_{0 j_m}(y^m)\right]\neq\emptyset$ for $y, y^m,m=1,\ldots,q-1\in\R^3$. There is a cube, denoted by $\Q\subset \R^2$, such that
\bel{inclusion}
\Q~\subset~\Lambda_{r j}^*(y_1,y_2)\cap \left[\Cap_{m=1}^{q-1} \Lambda_{0 j_m}^*(y^m_1,y^m_2)\right]
\eeq
and
\bel{cube size}
\vol\left\{ \Q\right\}~\approx~2^{j_\nu}2^{j_\nu}.
\eeq
\end{lemma}

{\bf Proof}: Observe that $\Lambda_{r j}^1(y_1)=\Lambda_{\ell j}^1(y_1)$. 
Because $\Gamma_{\ell j}(y)\cap \left[\Cap_{m=1}^{q-1} \Gamma_{0 j_m}(y^m)\right]$ is non-empty, there is an 
$(\Hat{x}_1,\Hat{x}_2)\in \Lambda_{\ell j}(y_1,y_2)\cap\left[\Cap_{m=1}^{q-1}\Lambda_{0 j_m}(y^m_1,y^m_2)\right]$ such that
\bel{dist est.1}
\begin{array}{cc}\ds
\Hat{x}_1~\in~\Lambda_{r j}^1(y_1)\cap\left[\Cap_{m=1}^{q-1}\Lambda_{0 j_m}^1(y^m_1)\right],
\\\\ \ds
\left|\Hat{x}_2-y_2\right|~<~2^{j-\ell+1}~\leq~ 2^{j_\nu-2},\qquad 2^{j_m}~\leq~ \left|\Hat{x}_2-y^m_2\right|~<~2^{j_m+1},\qquad m=1,2,\ldots,q-1.
\end{array}
\eeq
By using (\ref{dist est.1}) and the triangle inequality, we find
\bel{dist est.2}
\begin{array}{lr}\ds
2^{j_m-1}~<~2^{j_m}-2^{j_\nu-2}~<~\left|\Hat{x}_2-y^m_2\right|-\left| \Hat{x}_2-y_2\right|~\leq~\left| y_2-y^m\right|,
\\\\ \ds
\left| y_2-y^m\right|~\leq~\left|\Hat{x}_2-y^m_2\right|+\left| \Hat{x}_2-y_2\right|~<~2^{j_m+1}+2^{j_\nu-2}~<~2^{j_m+2}.
\end{array}
\eeq
For any $x_2\in\Lambda_{0 j_\nu-3}^2(y_2)$, we have $|x_2-y_2|<2^{j_\nu-2}$. By using (\ref{dist est.2}) and the triangle inequality again, we find
\bel{dist est.3}
\begin{array}{lr}\ds
2^{j_m-3}~<~2^{j_m-1}-2^{j_\nu-2}~<~\left|y_2-y^m_2\right|-|x_2-y_2|~\leq~\left|x_2-y_2^m\right|,
\\\\ \ds
\left|x_2-y_2^m\right|~\leq~\left|y_2-y^m_2\right|+|x_2-y_2|~<~2^{j_m+2}+2^{j_\nu-2}~<~2^{j_m+3}.
\end{array}
\eeq
This implies $x_2\in {^*}\Lambda_{0 j_m}^2(y^m_2)$ for every $m=1,\ldots,q-1$. Recall $r=j-j_\nu+2$. We have $\Lambda_{0 j_\nu-2}^2(y_2)=\Lambda_{r j}^2(y_2)$. From the above estimates, we obtain
\bel{Q_2 inclusion}
\Lambda_{0 j_\nu-3}^2(y_2)~\subset~ {^*}\Lambda_{r j}^2(y_2)\cap\left[\Cap_{m=1}^{q-1}{^*}\Lambda_{0 j_m}^2(y^m_2)\right].
\eeq
Let $\Q_1\subset\Lambda_{0\j_\nu}^1(y^\nu_1)$ be an interval containing $\Hat{x}_1$ whose side length equals $2^{j_\nu-3}$. Clearly, $\Q_1$ intersects $\Lambda_{r j}^1(y_1)$ and every $\Lambda_{0 j_m}^1(y^m_1), m=1,\ldots,q-1$. 
Consequently, we must have 
\bel{Q_1 inclusion}
\Q_1~\subset~ {^*}\Lambda_{r j}^1(y_1)\cap\left[\Cap_{m=1}^{q-1}{^*}\Lambda_{0 j_m}^1(y^m_1)\right].
\eeq
Define $\Q=\Q_1\times \Lambda_{0 j_\nu-3}^2(y_2)$. From (\ref{Q_2 inclusion})-(\ref{Q_1 inclusion}), we conclude (\ref{inclusion})-(\ref{cube size}).\endproof

Let $\Lambda_{\ell j}(x_1,x_2)$ defined in (\ref{Gamma project}). Consider
\bel{KERNEL Sum}
\begin{array}{lr}\ds
\int_{\Gamma_{\ell j}(y)\cap \left[\Cap_{m=1}^{q-1} \Gamma_{0 j_m}(y^m)\right]}    \V^{\alpha\beta\vartheta}(x-y)\prod_{m=1}^{q-1} \V^{\alpha\beta}(x-y^m)dx
\\\\ \ds
~=~\sum_{k\in\Z} \int_{2^{j}2^{j-\ell}2^{-k}\leq|x_3-y_3|<2^{j}2^{j-\ell}2^{-k+1}}
\\ \ds~~~~~~~~~~~~~~
\left\{\iint_{\Lambda_{\ell j}(y_1,y_2)\cap \left[\Cap_{m=1}^{q-1} \Lambda_{0 j_m}(y^m_1, y^m_2)\right]}    \V^{\alpha\beta\vartheta}(x-y)\prod_{m=1}^{q-1} \V^{\alpha\beta}(x-y^m)dx_1 dx_2\right\}dx_3.
\end{array}
\eeq
By definition of $\Lambda_{\ell j}(x_1,x_2)$, we find
\bel{projection volume}
\vol\Bigg\{ \Lambda_{\ell j}(y_1,y_2)\cap \left[\Cap_{m=1}^{q-1} \Lambda_{0 j_m}(y^m_1, y^m_2)\right]\Bigg\}~\lesssim~2^{j-\ell} 2^{j_\nu}.
\eeq
Recall $\Gamma_{\ell j k}(x)$ defined in (\ref{Gamma_ljk}) for $\ell, j,k\in\Z$. We have
\bel{KERNEL k}
\begin{array}{lr}\ds
\int_{2^{j}2^{j-\ell}2^{-k}\leq|x_3-y_3|<2^{j}2^{j-\ell}2^{-k+1}}

\left\{\iint_{\Lambda_{\ell j}(y_1,y_2)\cap \left[\Cap_{m=1}^{q-1} \Lambda_{0 j_m}(y^m_1, y^m_2)\right]}    \V^{\alpha\beta\vartheta}(x-y)\prod_{m=1}^{q-1} \V^{\alpha\beta}(x-y^m)dx_1 dx_2\right\}dx_3
\\\\ \ds
~=~ \int_{\Gamma_{\ell j k}(y)\cap \left[\Cap_{m=1}^{q-1} \Gamma_{0 j_mk_m}(y^m)\right]}    \V^{\alpha\beta\vartheta}(x-y)\prod_{m=1}^{q-1} \V^{\alpha\beta}(x-y^m)dx
\end{array}
\eeq
where 
\bel{k formula}
k_m~=~k-\big[ (j-j_m)+(j-j_m)-\ell\big],\qquad m~=~1,2,\ldots,q-1.
\eeq
\begin{remark}
Note that $\Gamma_{\ell j k}(y)= \Lambda_{\ell j}(y_1,y_2)\times \left\{ x_3\in\R\colon 2^{j}2^{j-\ell}2^{-k}\leq|x_3-y_3|<2^{j}2^{j-\ell}2^{-k+1}\right\}$ and
$\Gamma_{0 j_m k_m}(y^m)= \Lambda_{0 j_m}(y_1^m,y_2^m)\times \left\{ x_3\in\R\colon 2^{j_m}2^{j_m}2^{-k_m}\leq|x_3-y_3|<2^{j_m}2^{j_m}2^{-k_m+1}\right\}$. We thus have 
\[2^{j}2^{j-\ell}2^{-k}~=~2^{j_m}2^{j_m}2^{-k_m},\qquad m~=~1,2,\ldots,q-1.\]
\end{remark}
Let $\V^{\alpha\beta\vartheta}(x)$ defined in (\ref{V}). We have
\bel{V EST}
\begin{array}{lr}\ds
\V^{\alpha\beta\vartheta}(x-y)\prod_{m=1}^{q-1} \V^{\alpha\beta}(x-y^m)~\approx~
\\\\ \ds
\Big[2^j 2^{j-\ell}\Big]^{\alpha+\beta-2} 2^{-k(\beta-1)} \left[2^k+2^{-k}\right]^{-\vartheta}\prod_{m=1}^{q-1} \Big[2^{j_m} 2^{j_m}\Big]^{\alpha+\beta-2}2^{-k_m(\beta-1)}  \left[2^{k_m}+2^{-k_m}\right]^{-1} 
\end{array}
\eeq
whenever $x\in \Gamma_{\ell j k}(y)\cap \left[\Cap_{m=1}^{q-1} \Gamma_{0 j_m k_m}(y^m)\right]$.

By using (\ref{projection volume}) and (\ref{V EST}), we have
\bel{KERNEL Est.1}
\begin{array}{lr}\ds
\int_{\Gamma_{\ell j k}(y)\cap \left[\Cap_{m=1}^{q-1} \Gamma_{0 j_mk_m}(y^m)\right]}    \V^{\alpha\beta\vartheta}(x-y)\prod_{m=1}^{q-1} \V^{\alpha\beta}(x-y^m)dx
\\\\ \ds
~\approx~\Big[2^j 2^{j-\ell}\Big]^{\alpha+\beta-2} 2^{-k(\beta-1)}\left[2^k+2^{-k}\right]^{-\vartheta} \prod_{m=1}^{q-1} \Big[2^{j_m} 2^{j_m}\Big]^{\alpha+\beta-2}2^{-k_m(\beta-1)}  \left[2^{k_m}+2^{-k_m}\right]^{-1}  
\Big[2^{j_\nu} 2^{j-\ell}\Big] 2^j 2^{j-\ell} 2^{-k}
\\\\ \ds
~=~\Big[2^j 2^{j-\ell}\Big]^{\alpha+\beta} 2^{-k(\beta-1)} \left[2^k+2^{-k}\right]^{-\vartheta}\prod_{m=1}^{q-1} \Big[2^{j_m} 2^{j_m}\Big]^{\alpha+\beta-2}2^{-k_m(\beta-1)}  \left[2^{k_m}+2^{-k_m}\right]^{-1}  
\Big[2^{j_\nu-j} \Big]  2^{-k}
\\\\ \ds
~=~2^{(r-\ell)(\alpha+\beta)}\Big[2^j 2^{j-r}\Big]^{\alpha+\beta} 2^{-k(\beta-1)} \left[2^k+2^{-k}\right]^{-\vartheta}\prod_{m=1}^{q-1} \Big[2^{j_m} 2^{j_m}\Big]^{\alpha+\beta-2}2^{-k_m(\beta-1)}  \left[2^{k_m}+2^{-k_m}\right]^{-1}  
\Big[2^{j_\nu-j} \Big]  2^{-k}
\\\\ \ds
~=~2^{(r-\ell)(\alpha+\beta)}\Big[2^j 2^{j-r}\Big]^{\alpha+\beta-2} 2^{-k(\beta-1)} \left[2^k+2^{-k}\right]^{-\vartheta}\prod_{m=1}^{q-1} \Big[2^{j_m} 2^{j_m}\Big]^{\alpha+\beta-2}2^{-k_m(\beta-1)}  \left[2^{k_m}+2^{-k_m}\right]^{-1}  
\\ \ds~~~~~~~
\Big[2^{j_\nu} 2^{j-r}\Big] 2^j 2^{j-r} 2^{-k}
\\\\ \ds
~=~2^{(r-\ell)(\alpha+\beta)}\Big[2^j 2^{j-r}\Big]^{\alpha+\beta-2} 2^{-k\beta} \left[2^k+2^{-k}\right]^{-\vartheta}\prod_{m=1}^{q-1} \Big[2^{j_m} 2^{j_m}\Big]^{\alpha+\beta-2}2^{-k_m(\beta-1)}  \left[2^{k_m}+2^{-k_m}\right]^{-1}  
\\ \ds~~~~~~~
\Big[2^{j_\nu} 2^{j_\nu-2}\Big] 2^j 2^{j-r}.\qquad \hbox{\small{( $r=j-j_\nu-2$ )}}
\end{array}
\eeq
On the other hand, recall ${^*}\Gamma_{\ell j}(x)$ and ${^*}\Lambda_{\ell j}(x_1,x_2)$ defined in (\ref{Gamma*_ljk})-(\ref{Gamma*_lj}) and (\ref{Gamma_ljk*}).
We have
\bel{KERNEL Est.2}
\begin{array}{lr}\ds
\int_{{^*}\Gamma_{r j}(y)\cap \left[\Cap_{m=1}^{q-1} {^*}\Gamma_{0 j_m}(y^m)\right]}    \V^{\alpha\beta\vartheta}(x-y)\prod_{m=1}^{q-1} \V^{\alpha\beta}(x-y^m)dx
\\\\ \ds
~=~\sum_{k\in\Z} \int_{2^{j}2^{j-\ell}2^{-k}\leq|x_3-y_3|<2^{j}2^{j-\ell}2^{-k+1}}
\\ \ds~~~~~~~~~~~~~
\left\{\iint_{{^*}\Lambda_{r j}(y_1,y_2)\cap \left[\Cap_{m=1}^{q-1} {^*}\Lambda_{0 j_m}(y^m_1, y^m_2)\right]}    \V^{\alpha\beta\vartheta}(x-y)\prod_{m=1}^{q-1} \V^{\alpha\beta}(x-y^m)dx_1 dx_2\right\}dx_3.
\end{array}
\eeq
Let $\kappa=k-(r-\ell)$. Each integral in the summand of (\ref{KERNEL Est.2}) can be written as
\bel{Case3 Est k to kappa}
\begin{array}{lr}\ds
\int_{2^{j}2^{j-r}2^{-\kappa}\leq|x_3-y_3|<2^{j}2^{j-r}2^{-\kappa+1}}\left\{\iint_{{^*}\Lambda_{r j}(y_1,y_2)\cap \left[\Cap_{m=1}^{q-1} {^*}\Lambda_{0 j_m}(y^m_1, y^m_2)\right]}    \V^{\alpha\beta\vartheta}(x-y)\prod_{m=1}^{q-1} \V^{\alpha\beta}(x-y^m)dx_1 dx_2\right\}dx_3
\\\\ \ds
~=~\int_{{^*}\Gamma_{r j \kappa}(y)\cap \left[\Cap_{m=1}^{q-1} {^*}\Gamma_{0 j_m\kappa_m}(y^m)\right]}    \V^{\alpha\beta\vartheta}(x-y)\prod_{m=1}^{q-1} \V^{\alpha\beta}(x-y^m)dx
\end{array}
\eeq
where 
\bel{kappa_m}
\begin{array}{lr}\ds
\kappa_m~=~\kappa-\Big[ (j-j_m)+(j-j_m)-r\Big]
\\\\ \ds~~~~~
~=~k-\Big[ (j-j_m)+(j-j_m)-\ell\Big]~=~k_m, \qquad m=1,2,\ldots,q-1\qquad \hbox{\small{by (\ref{k formula})}}.
\end{array}
\eeq
Recall  {\bf Lemma 5.1}. By using (\ref{V EST}), 
we have
\bel{KERNEL Est.3}
\begin{array}{lr}\ds
\int_{{^*}\Gamma_{r j \kappa}(y)\cap \left[\Cap_{m=1}^{q-1} {^*}\Gamma_{0 j_m \kappa_m}(y^m)\right]}    \V^{\alpha\beta\vartheta}(x-y)\prod_{m=1}^{q-1} \V^{\alpha\beta}(x-y^m)dx
\\\\ \ds
~\approx~\Big[2^j 2^{j-r}\Big]^{\alpha+\beta-2} 2^{-\kappa(\beta-1)} \left[2^\kappa+2^{-\kappa}\right]^{-\vartheta}\prod_{m=1}^{q-1} \Big[2^{j_m} 2^{j_m}\Big]^{\alpha+\beta-2}2^{-\kappa_m(\beta-1)}  \left[2^{\kappa_m}+2^{-\kappa_m}\right]^{-1} 
\\\\ \ds~~~~~~~~~~~~~~~~~~~~~~~~~~~~~~~~~~~~~~~~~~
\int_{2^j 2^{j-r}2^{-\kappa}}^{2^j 2^{j-r}2^{-\kappa+1}}\left\{\iint_{{^*}\Lambda_{r j}(y_1,y_2)\cap \left[\Cap_{m=1}^{q-1} {^*}\Lambda_{0 j_m}(y^m_1,y^m_2)\right]}    dx_1 dx_2\right\} dx_3
\\\\ \ds
~\ge~
\Big[2^j 2^{j-r}\Big]^{\alpha+\beta-2} 2^{-\kappa(\beta-1)} \left[2^\kappa+2^{-\kappa}\right]^{-\vartheta}\prod_{m=1}^{q-1} \Big[2^{j_m} 2^{j_m}\Big]^{\alpha+\beta-2}2^{-\kappa_m(\beta-1)}  \left[2^{\kappa_m}+2^{-\kappa_m}\right]^{-1}  
\\ \ds~~~~~~~~~~~~~~~~~~~~~~~~~~~~~~~~~~~~~~~~~~~
\int_{2^j 2^{j-r}2^{-\kappa}}^{2^j 2^{j-r}2^{-\kappa+1}}\left\{ \iint_{\Q}    dx_1 dx_2\right\}dx_3     \qquad\hbox{\small{by (\ref{inclusion})}}
  \\\\ \ds
 ~\approx~\Big[2^j 2^{j-r}\Big]^{\alpha+\beta-2} 2^{-\kappa(\beta-1)} \left[2^\kappa+2^{-\kappa}\right]^{-\vartheta}\prod_{m=1}^{q-1} \Big[2^{j_m} 2^{j_m}\Big]^{\alpha+\beta-2}2^{-\kappa_m(\beta-1)}  \left[2^{\kappa_m}+2^{-\kappa_m}\right]^{-1} 
 \\ \ds~~~~~~~~
 2^j 2^{j-r} 2^{-\kappa}2^{j_\nu}2^{j_\nu} \qquad\hbox{\small{by (\ref{cube size})}}
 \\\\ \ds
 ~=~\Big[2^j 2^{j-r}\Big]^{\alpha+\beta-2} 2^{-k\beta}2^{(r-\ell)\beta} \left[2^{k-(r-\ell)}+2^{-[k-(r-\ell)]}\right]^{-\vartheta}\prod_{m=1}^{q-1} \Big[2^{j_m} 2^{j_m}\Big]^{\alpha+\beta-2}2^{-k_m(\beta-1)}  \left[2^{k_m}+2^{-k_m}\right]^{-1} 
 \\ \ds~~~~~~~~
 2^j 2^{j-r} 2^{j_\nu}2^{j_\nu} 
 \\\\ \ds
 ~\gtrsim~2^{(r-\ell)(\beta+\vartheta)}\Big[2^j 2^{j-r}\Big]^{\alpha+\beta-2} 2^{-k\beta} \left[2^{k}+2^{-k}\right]^{-\vartheta}
 \\ \ds~~~~~~~
 \prod_{m=1}^{q-1} \Big[2^{j_m} 2^{j_m}\Big]^{\alpha+\beta-2}2^{-k_m(\beta-1)}  \left[2^{k_m}+2^{-k_m}\right]^{-1} 
 2^j 2^{j-r} 2^{j_\nu}2^{j_\nu}.  
\end{array}
\eeq
From (\ref{KERNEL Est.1}) and (\ref{KERNEL Est.3}), we obtain
\bel{KERNELS <}
\begin{array}{lr}\ds
\int_{\Gamma_{\ell j}(y)\cap \left[\Cap_{m=1}^{q-1} \Gamma_{0 j_m}(y^m)\right]}    \V^{\alpha\beta\vartheta}(x-y)\prod_{m=1}^{q-1} \V^{\alpha\beta}(x-y^m)dx
\\\\ \ds
~\leq~\C_\beta~2^{(r-\ell)(\alpha-\vartheta)}~\int_{{^*}\Gamma_{r j}(y)\cap \left[\Cap_{m=1}^{q-1} {^*}\Gamma_{0 j_m}(y^m)\right]}    \V^{\alpha\beta\vartheta}(x-y)\prod_{m=1}^{q-1} \V^{\alpha\beta}(x-y^m)dx
\end{array}
\eeq
where $\alpha-\vartheta>0$.

Recall (\ref{Int change}). By using (\ref{KERNELS <}), we have
\bel{Int change EST}
\begin{array}{lr}\ds
\int_{\R^3} \left\{\int_{\Gamma_{\ell j}(x)} f(y)\V^{\alpha\beta\vartheta}(x-y)dy\right\}\prod_{m=1}^{q-1} \left\{\int_{\Gamma_{0 j_m}(x)} f(y^m)\V^{\alpha\beta}(x-y^m)dy^m\right\} dx~=~
\\\\ \ds
\idotsint_{\R^3\times\cdots\times\R^3}f(y)\prod_{m=1}^{q-1} f(y^m) \left\{\int_{\Gamma_{\ell j}(y)\cap \left[\Cap_{m=1}^{q-1} \Gamma_{0 j_m}(y^m)\right]}    \V^{\alpha\beta\vartheta}(x-y)\prod_{m=1}^{q-1} \V^{\alpha\beta}(x-y^m)dx\right\} dy \prod_{m=1}^{q-1} dy^m
\\\\ \ds
~\leq~\C_\beta~2^{(r-\ell)(\alpha-\vartheta)}~\idotsint_{\R^3\times\cdots\times\R^3}f(y)\prod_{m=1}^{q-1} f(y^m)
\\ \ds~~~~~~~~~~~~~~~~~~~~~~~~~~~~~~
\left\{ \int_{{^*}\Gamma_{r j}(y)\cap \left[\Cap_{m=1}^{q-1} {^*}\Gamma_{0 j_m}(y^m)\right]}    \V^{\alpha\beta\vartheta}(x-y)\prod_{m=1}^{q-1} \V^{\alpha\beta}(x-y^m)dx\right\} dy \prod_{m=1}^{q-1} dy^m
\\\\ \ds
~=~\C_\beta~2^{(r-\ell)(\alpha-\vartheta)}\int_{\R^3} \left\{\int_{{^*}\Gamma_{r j}(x)} f(y)\V^{\alpha\beta\vartheta}(x-y)dy\right\}\prod_{m=1}^{q-1} \left\{\int_{{^*}\Gamma_{0 j_m}(x)} f(y^m)\V^{\alpha\beta}(x-y^m)dy^m\right\} dx.
\end{array}
\eeq
Note that $r=j-j_\nu+2>0$ because $j>j_\nu$. We find $j-r=j_\nu-2$. This brings us back to the situation of  {\bf Case 1}. 
By carrying out the estimates in analogue to (\ref{theta Case1 est.1})-(\ref{Case1.Est.2}), we find
\bel{Case 3 Est}
\begin{array}{lr}\ds
\int_{\R^3} \left\{\int_{\Gamma_{\ell j}(x)} f(y)\V^{\alpha\beta\vartheta}(x-y)dy\right\}\prod_{m=1}^{q-1} \left\{\int_{\Gamma_{0 j_m}(x)} f(y^m)\V^{\alpha\beta}(x-y^m)dy^m\right\} dx
\\\\ \ds
~\leq~\C_\beta~2^{(r-\ell)(\alpha-\vartheta)}\int_{\R^3} \left\{\int_{{^*}\Gamma_{r j}(x)} f(y)\V^{\alpha\beta\vartheta}(x-y)dy\right\}\prod_{m=1}^{q-1} \left\{\int_{{^*}\Gamma_{0 j_m}(x)} f(y^m)\V^{\alpha\beta}(x-y^m)dy^m\right\} dx
\\\\ \ds
~\leq~\C_{\alpha~\beta~q}~2^{(r-\ell)(\alpha-\vartheta)} ~2^{-{1\over 3}\sigma r}
~
2^{-\left| j-\lambda(r,x)\right| \sigma} 
\prod_{m=1}^{q-1} 2^{-\left|j_m-\lambda(0,x)\right| \sigma} 
\\\\ \ds~~~~~~~~~~~~~~
\left\|f\right\|_{\L^p(\R^3)}^{q-p}
\int_{\R^3}\Big(\varphi_0(x)\Big)^{(q-2)\big[{1\over p}-{1\over q}\big]} \Big(\M f\Big)^p(x)dx.
\end{array}
\eeq
Recall $\sigma=\min\big\{\alpha+\beta, {2\over q}\big\}$. Choose $\vex=\min\{ \alpha-\vartheta, \sigma\}$.

From (\ref{Case 3 Est}), we further have
\bel{Case3 EST}
\begin{array}{lr}\ds
\sum_{\G_3}\int_{\R^3} \left\{\int_{\Gamma_{\ell j}(x)} f(y)\V^{\alpha\beta\vartheta}(x-y)dy\right\}\prod_{m=1}^{q-1} \left\{\int_{\Gamma_{0 j_m}(x)} f(y^m)\V^{\alpha\beta}(x-y^m)dy^m\right\} dx
\\\\ \ds
~\leq~\C_{\alpha~\beta~q}~2^{{1\over 3}(r-\ell)\vex} ~2^{-{1\over 3}\vex r} \sum_{j,j_1,\ldots,j_{q-1}\in\Z} 2^{-\left| j-\lambda(r,x)\right| \vex} 
\prod_{m=1}^{q-1} 2^{-\left|j_m-\lambda(0,x)\right| \vex} 
\\ \ds~~~~~~~~~~~~~~~~~~~~~~~~~~~~~~~~~~~~~~~~~~~~~~~~~~~~~
\left\|f\right\|_{\L^p(\R^3)}^{q-p}
\int_{\R^3}\Big(\varphi_0(x)\Big)^{(q-2)\big[{1\over p}-{1\over q}\big]} \Big(\M f\Big)^p(x)dx
\\\\ \ds
~\leq~\C_{\alpha~\beta~q}~2^{-{1\over 3}\vex\ell} ~\left\|f\right\|_{\L^p(\R^3)}^{q-p}
\int_{\R^3}\Big(\varphi_0(x)\Big)^{(q-2)\big[{1\over p}-{1\over q}\big]} \Big(\M f\Big)^p(x)dx
\end{array}
\eeq
where $\G_3=\left\{ j,j_1,\ldots,j_{q-1}\in\Z\colon j-\ell< j_\nu-2<j-2\right\}$.

{\small Department of Mathematics, Westlake University}~~~~~~~~~~~~~~~
{\small wangzipeng@westlake.edu.cn}

\end{document}